\documentclass[12pt,a4paper]{article}
\usepackage[]{fontenc}
\usepackage[latin1]{inputenc}
\usepackage[francais]{babel}
\usepackage{a4,amsmath,amssymb,amsfonts,amsthm}
\usepackage{amsrefs}
\usepackage[all]{xy}

\begin{document}
\bibliographystyle{plain}
 

\def\mR{\M{R}}           
\def\mZ{\M{Z}}           
\def\mN{\M{N}}           
\def\mQ{\M{Q}}       
\def\mC{\M{C}}  
\def\mG{\M{G}}



\def\Spec{{\rm Spec}}
\def\rg{{\rm rg}}
\def\Hom{{\rm Hom}}
\def\Aut{{\rm Aut}}
 \def\Tr{{\rm Tr}}
 \def\Exp{{\rm Exp}}
 \def\Gal{{\rm Gal}}
 \def\End{{\rm End}}
 \def\det{{{\rm det}}}
 \def\Td{{\rm Td}}
 \def\ch{{\rm ch}}
 \def\che{{\rm ch}_{\rm eq}}
  \def\Spec{{\rm Spec}}
\def\Id{{\rm Id}}
\def\Zar{{\rm Zar}}
\def\Supp{{\rm Supp}}
\def\eq{{\rm eq}}
\def\Ann{{\rm Ann}}
\def\LT{{\rm LT}}
\def\Pic{{\rm Pic}}
\def\rg{{\rm rg}}
\def\et{{\rm et}}
\def\sep{{\rm sep}}
\def\ppcm{{\rm ppcm}}
\def\ord{{\rm ord}}


\def\beginProof{\par{\bf Preuve. }}
 \def\endProof{${\qed}$\par\smallskip}
 \def\pr{^{\prime}}
 \def\prpr{^{\prime\prime}}
 \def\mtr#1{\overline{#1}}
 \def\ra{\rightarrow}
 \def\mfp{{\mathfrak p}}
 
 \def\mQ{{\Bbb Q}}
 \def\mR{{\Bbb R}}
 \def\mZ{{\Bbb Z}}
 \def\mC{{\Bbb C}}
 \def\mN{{\Bbb N}}
 \def\mF{{\Bbb F}}
 \def\mA{{\Bbb A}}
  \def\mG{{\Bbb G}}
 \def\CI{{\cal I}}
 \def\CH{{\cal H}}
 \def\CO{{\cal O}}
 \def\CA{{\cal A}}
 \def\CB{{\cal B}}
 \def\CC{{\cal C}}
 \def\CK{{\cal K}}
 \def\CL{{\cal L}}
 \def\CM{{\cal M}}
\def\CP{{\cal P}}
 \def\wt#1{\widetilde{#1}}
 \def\mod{{\rm mod\ }}
 \def\refeq#1{(\ref{#1})}
 \def\blb{{\big(}}
 \def\brb{{\big)}}
\def\mc{{{\mathfrak c}}}
\def\mcpr{{{\mathfrak c}'}}
\def\mcprpr{{{\mathfrak c}''}}
\def\ul#1{\overline{#1}}
\def\ss{{\rm ss}}
\def\parf{{\rm parf}}
\def\P1{{{\bf P}^1}}
\def\cod{{\rm cod}}
\def\pr{\prime}
\def\prpr{\prime\prime}
\def\ss{\scriptstyle}
\def\OX{{ {\cal O}_X}}
\def\mpartial{{\mtr{\partial}}}
\def\inv{{\rm inv}}
\def\indlim{\underrightarrow{\lim}}
\def\prolim{\underleftarrow{\lim}}
\def\pprolim{'\prolim'}
\def\Pro{{\rm Pro}}
\def\Ind{{\rm Ind}}
\def\Ens{{\rm Ens}}
\def\without{\backslash}
\def\pbdb{{\Pro_b\ D^-_c}}
\def\qc{{\rm qc}}
\def\Com{{\rm Com}}
\def\an{{\rm an}}
\def\gfield{{\rm\bf k}}
\def\s{{\rm s}}
\def\dR{{\rm dR}}
\def\ari#1{\widehat{#1}}
\def\ul#1{\underbar{#1}}
\def\sul#1{\underline{\scriptsize #1}}
\def\mou{{\mathfrak u}}
\def\ich{\mathfrak{ch}}
\def\cl{{\rm cl}}
\def\K{{\rm K}}
\def\R{{\rm R}}
\def\F{{\rm F}}
\def\L{{\rm L}}
\def\pgcd{{\rm pgcd}}
\def\rc{{\rm c}}
\def\N{{\rm N}}
\def\E{{\rm E}}
\def\H{{\rm H}}
\def\CHOW{{\rm CH}}
\def\A{{\rm A}}
\def\d{{\rm d}}
\def\Res{{\rm  Res}}
\def\GL{{\rm GL}}
\def\Alb{{\rm Alb}}
\def\alb{{\rm alb}}
\def\Hdg{{\rm Hdg}}
\def\Num{{\rm Num}}
\def\Irr{{\rm Irr}}
\def\Frac{{\rm Frac}}


\def\RHom{{\rm RHom}}
\def\rRHom{{\mathcal RHom}}
\def\rHom{{\mathcal Hom}}
\def\dotimes{{\overline{\otimes}}}
\def\Ext{{\rm Ext}}
\def\rExt{{\mathcal Ext}}
\def\Tor{{\rm Tor}}
\def\rTor{{\mathcal Tor}}
\def\SP{{\mathfrak S}}

\def\H{{\rm H}}
\def\D{{\rm D}}
\def\Del{{\mathfrak D}}

 \newtheorem{theor}{Théorème}[section]
 \newtheorem{prop}[theor]{Proposition}
 \newtheorem{propdef}[theor]{Proposition-Définition}
 \newtheorem{cor}[theor]{Corollaire}
 \newtheorem{lemme}[theor]{Lemme}
 \newtheorem{slemme}[theor]{sous-lemme}
 \newtheorem{defin}[theor]{Définition}
 \newtheorem{conj}[theor]{Conjecture}

 \parindent=0pt
 \parskip=5pt

 \author{
 Vincent MAILLOT\footnote{Institut de Math\'ematiques de Jussieu,
 Universit\'e Paris 7 Denis Diderot, C.N.R.S.,
 Case Postale 7012,
 2 place Jussieu,
 F-75251 Paris Cedex 05, France,
 E-mail : vmaillot@math.jussieu.fr}\ \ et 
 Damian RÖSSLER\footnote{\'Equipe d'Arithmétique et Géométrie Algébrique, 
Département de Mathématiques, 
Bâtiment 425, 
Faculté des Sciences d'Orsay, 
Université Paris-Sud 11, 
F-91405 Orsay Cedex, 
France, E-mail: damian.rossler@math.u-psud.fr}}
 \title{Une conjecture sur la torsion des classes de Chern des fibrés de Gauss-Manin}
 \date{}
\maketitle
\begin{abstract}
Pour tout $t\in\mN$ nous définissons un certain entier positif $\N_t$ 
et nous conjecturons: si $H$ est un fibré de Gauss-Manin d'une fibration semi-stable alors 
la $t$-ème 
classe de Chern de $H$  
est annulée par $\N_t$. 
Nous démontrons diverses conséquences de cette conjecture.  
\end{abstract}

 \section{Introduction}

 Soit $f:X\ra Y$ un morphisme projectif de 
 schémas lisses et quasi-projectifs sur un corps algébriquement clos $K$ de caractéristique nulle. Soit $D\hookrightarrow X$ 
 et $E\hookrightarrow Y$ des diviseurs à croisements normaux. 
 On suppose que $f$ est semi-stable relativement à $D$ et $E$. 
 Voir \cite[Par. 1, Def. 1.1]{Illusie-Reduction} ou la sous-section \ref{morss} plus bas pour la définition d'un diviseur à  croisements normaux et d'un morphisme semi-stable.
 
 On notera $\Omega^\bullet_{X/Y}(\log)$ 
 le complexe de de Rham logarithmique relatif de $X$ au-dessus de $Y$ relativement 
 à $D$ et $E$. Voir  \cite[Par. 1, Def. 1.3]{Illusie-Reduction} ou la sous-section \ref{morss} plus bas pour la définition de cette notion. 
 
 Pour tout $j\geqslant 0$, on écrira  
 $$\H_\dR^{j}(X/Y)(\log):=\R^{j} f_*(\Omega^\bullet_{X/Y}(\log))$$ pour le $j$-ème 
 fibré de Gauss-Manin logarithmique relativement à $f$, $D$ et $E$. 
 On peut montrer que $\H_\dR^{j}(X/Y)(\log)$ est un faisceau localement 
 libre de rang fini. 
 
Enfin, pour tout $t\in\mN$, on définit 
\begin{displaymath}
\N_t:=
\left\{
\begin{array}{ll|}
0 & \textrm{si $t=0$}\\
2\cdot\prod_{p-1|t}p^{\ord_p(t)+1} & \textrm{si $t$ est pair et positif}\\
2 & \textrm{si $t$ est impair et positif}
\end{array}
\right.
\end{displaymath}
où $p$ parcourt les nombres premiers. 
Le théorème de von Staudt montre que si $t$ est pair et positif, alors  $\N_t/2={\rm D\acute{e}nominateur}((-1)^{t+2\over 2}B_t/t)$, 
où $B_t$ est le $t$-ème nombre de Bernoulli 
(cf. \cite[Appendix B]{Milnor-Stasheff}). On rappelle 
que le $t$-ème nombre 
de Bernoulli est défini par l'identité de séries formelles
$$
{t\over \exp(t)-1}=\sum_{n\geqslant 0}B_n {t^n\over n!}
$$
L'objet du présent 
 article est de proposer la conjecture suivante:
 \begin{conj} Pour tout $j,t\geqslant 0$,  l'égalité
 $$
 \N_t\cdot \rc_{t}(\H_\dR^{j}(X/Y)(\log))=0
 $$
 est vérifiée dans ${\rm CH}^t(Y)$. 
 \label{mainconj}
 \end{conj}
Ici, $\CHOW^t(Y)$ est le $t$-ième groupe de Chow de $Y$ 
(cf. \cite{Fulton-Intersection}) et $\rc_t(\bullet)$ est la $t$-ème classe de Chern 
(à valeurs dans $\CHOW^t(Y)$). 

Les conséquences de la Conjecture \ref{mainconj} déjà démontrées dans la littérature 
mathématique sont les suivantes. 
Dans la situation où $f$ est lisse, 
il est démontré dans \cite{Grothendieck-Classes} (voir les 
calculs faits dans \cite{Evens-Kahn-Chern}) que l'image 
de la Conjecture \ref{mainconj} par l'application classe de cycle 
est vérifiée (pour un généralisation au cas non-lisse, voir 
le point (a) du Théorème \ref{mainth} plus bas). 
Encore dans la situation où $f$ est lisse, le fait que 
$\rc_t(\H_\dR^{1}(X/Y))$ est de torsion pour tout $t\geqslant 0$ est démontré 
(pour la première fois ?) dans 
\cite{VDG-Cycles}.
Dans \cite{Esnault-Viehweg-Chern}, il est démontré que 
$\rc_t(\H_\dR^{1}(X/Y)(\log))$ est de torsion pour tout $t\geqslant 0$ et 
que l'image par l'application classe de cycle de 
$\rc_t(\H_\dR^{j}(X/Y)(\log))$ est de torsion pour tout $j,t\geqslant 0$.  
Lorsque $f$ est lisse, il est démontré dans \cite[Appendix]{Bismut-Eta}
 que les classes de Cheeger-Simons 
de $\H_\dR^{j}(X/Y)$ sont de torsion pour tout $j\geqslant 0$; si on suppose de plus que $Y$ est 
projectif sur $K$, ceci implique que les classes de Chern en cohomologie 
de Deligne de $\H_\dR^{j}(X/Y)$ sont de torsion pour tout $j\geqslant 0$; ce dernier 
énoncé est aussi une conséquence d'un théorème de Reznikov (cf. \cite{Reznikov-All}). 
Lorsque le morphisme est étale, les travaux de Fulton et MacPherson 
dans \cite[Cor. 19.3]{Fulton-MacPherson-Char} démontrent la Conjecture 
\ref{mainconj}. Dans \cite{Pappas-Integral}, Pappas démontre 
un théorème de Grothendieck-Riemann-Roch sans dénominateurs pour les morphismes projectifs et lisses. Si on applique ce théorème au complexe 
de de Rham de $f$, on obtient des énoncés d'annulation 
pour les classes de Chern des fibrés de Gauss-Manin. Cependant, 
ces énoncés dépendent a priori de la dimension relative de $f$. 

On notera que l'énoncé que $\rc_t(\H_\dR^{j}(X/Y)(\log))$ est 
de torsion pour tout $j,t\geqslant 0$ est déjà conjectural. Cette forme faible 
est implicite dans 
les travaux de Bloch, Esnault et Viehweg (par ex. \cite{Esnault-Viehweg-Chern} 
et \cite{Bloch-Esnault-Algebraic}). Elle 
est formellement proposée dans \cite[Par. 4.2, Conj. 3]{Maillot-Rossler-On-the}, 
dans le cas où $f$ est lisse. 
 
Du point de vue des auteurs, la forme faible de la conjecture est motivée par 
une conjecture en théorie d'Arakelov (cf. \cite[Conjecture 3.1]{Maillot-Rossler-Conjectures}) 
dont elle est une conséquence. Quant au nombre $\N_t$, 
il est motivé par le théorème de Grothendieck mentionné plus haut. 

La Conjecture \ref{mainconj} n'est pas optimale 
pour $j=0,1$. 

Cas $j=0$: dans \cite[Cor. 19.3]{Fulton-MacPherson-Char}, il est démontré que 
lorsque $f$ est un morphisme étale et que $t$ est pair, on a 
$$
{1\over 2}N_t\cdot\rc_t(H^0_\dR(X/Y))=0.
$$
En particulier, $12\cdot\rc_2(H^0_\dR(X/Y))=0.$

Cas $j=1$: si 
$f$ est lisse, la Conjecture 
\ref{mainconj} et le résultat principal de \cite{Maillot-Rossler-On-the} impliquent que 
\begin{equation}
(\prod_{p,p\leqslant t}p^{\ord_p(t)})
(\prod_{p,p>t}p^{\ord_p[\pgcd(\Num((2^t-1)B_t),N_t)]})\cdot 
\rc_t(H^1_\dR(X/Y))=0.
\label{j1eq}
\end{equation}
Ici $p$ parcourt les nombres premiers. On note $\ord_p(x)$ 
l'ordre du nombre premier $p$ dans le nombre 
naturel $x$ et $\Num(r)$ le numérateur du nombre rationnel 
$r$. Un cas particulier de l'équation \refeq{j1eq} est l'équation
\begin{equation}
8\cdot 
\rc_2(H^1_\dR(X/Y))=0.
\label{j1}
\end{equation}
Cependant $N_2=24>12>8$, d'où la non-optimalité lorsque $j=0,1$. 

Cela suggère qu'une conjecture optimale tiendrait compte 
du poids $j$  du fibré de Gauss-Manin.

Au sujet de la Conjecture \ref{mainconj}, 
nous démontrerons les résultats partiels suivants.

Si $K=\mC$, nous écrirons  
$$\cl:\CHOW^\bullet(Y)\to{\rm H}^{2\bullet}(Y(\mC),\mZ)$$
pour l'application qui associe à un cycle algébrique 
sa classe dans la  cohomologie singulière de $Y(\mC)$.  

Si $Y$ est projectif sur $K$, nous écrirons
$$
\alb:\CHOW^{d_Y}(Y)_0\to\Alb(Y)(K)
$$
pour l'application qui associe à un $0$-cycle algébrique 
de degré nul son image dans la variété d'Albanese de $Y$. 

Soit $l$ un nombre premier. 
Si $Y$ est projectif sur $K$, nous écrirons aussi 
$$
\lambda^t_l:\CHOW^{t}(Y)[l^\infty]\to \H^{2t-1}_\et(Y,\mQ_l/\mZ_l(t))
$$
pour l'application d'Abel-Jacobi de Bloch 
(voir \cite{Bloch-Torsion} pour la définition). Elle donne lieu à une application
$$
\lambda^t:\Tor(\CHOW^t(Y))\to\bigoplus_{l\ {\rm premier}}\H^{2t-1}_\et(Y,\mQ_l/\mZ_l(t)).
$$

\begin{theor}
Soit $d_Y:=\dim(Y)$ la dimension de $Y$. 
\begin{description}
\item[{\rm\bf (a)}] On suppose que $K=\mC$. Pour tout $j,t\geqslant 0$, l'égalité
 $$
 \cl\big[\ \N_t\cdot\rc_t(\H_\dR^{j}(X/Y)(\log))\ \big]=0
 $$
 est vérifiée dans ${\rm H}^{2t}(Y(\mC),\mZ)$. 
 \item[{\rm\bf (b)}]  On suppose que $Y$ est projectif sur $K$. 
 L'égalité
 $$
 \alb\big[\ \N_{d_Y}\cdot\rc_{d_Y}(\H_\dR^{j}(X/Y)(\log))\ \big]=0
 $$
 est alors vérifiée dans $\Alb(Y)(K)$ pour tout $j\geqslant 0$.  
 \item[{\rm\bf (c)}] On suppose que $\rc_{d_Y}(\H_\dR^{j}(X/Y)(\log))$ est 
 de torsion et que $Y$ est projectif sur $K$. Alors 
 $$\N_{d_Y}\cdot\rc_{d_Y}(\H_\dR^{j}(X/Y)(\log))=0.$$
\item[{\rm\bf (d)}] On suppose que les composantes irréductibles de 
$D$ et $E$ sont lisses sur $K$. Pour tout $t\geqslant 0$, l'égalité
 $$
 \N_t\cdot \rc_{t}\Big(\sum_{j\geqslant 0}(-1)^j\H_\dR^{j}(X/Y)(\log)\Big)=0
 $$
 est vérifiée dans ${\rm CH}^t(Y)$. 
  \item[{\rm\bf (e)}]  L'égalité
$$
 \N_t\cdot \rc_{t}(\H_\dR^{0}(X/Y)(\log))=0
 $$
est vérifiée pour tout $t\geqslant 0$. Si on suppose que les composantes irréductibles de 
$D$ et $E$ sont lisses sur $K$ et que les fibres  de $f$ sont de dimension $1$, alors
$$
 \N_t\cdot \rc_{t}(\H_\dR^{j}(X/Y)(\log))=0
 $$
 pour tout $j,t\geqslant 0$. 
 \item[${\rm\bf (f)}$] On suppose que $f$ est lisse, que 
 $K=\mC$ et que l'image de la représentation 
 du groupe fondamental de $Y(\mC)$ associée au système 
 localement constant $\R^j f(\mC)_*\mQ$ est finie. Alors
 $$
 \N_t\cdot \rc_{t}(\H_\dR^{j}(X/Y)(\log))=0
 $$
 pour tout $t\geqslant 0$. 
\item[${\rm\bf (g)}$] On suppose que $Y$ est projectif sur $K$. On suppose 
aussi que $\rc_{t}(\H_\dR^{j}(X/Y)(\log))$ est 
 de torsion. Alors
$$
 \lambda^t\big[\ \N_t\cdot\rc_t(\H_\dR^{j}(X/Y)(\log))\ \big]=0.
$$
\item[${\rm\bf (h)}$] On suppose que $Y$ est projectif sur $K$. 
On suppose 
aussi que $\rc_{2}(\H_\dR^{j}(X/Y)(\log))$ est 
 de torsion. Alors
 $$
 \N_2\cdot\rc_2(\H_\dR^{j}(X/Y)(\log))=0.
 $$
\end{description}
 \label{mainth}
 \end{theor} 

Les démonstrations des différents points du Théorème \ref{mainth} 
ont des structures semblables, même si les outils utilisés varient.
Dans les preuves de (a), (b), (d), (f) et (g), on démontre chaque fois que 
l'image de (d'une combinaison linéaire de) $\rc_t(\H^j_\dR(X/Y))$ 
dans un certain groupe est invariante 
par multiplication par $p^t$, pour presque tout $p$. 
Un lemme de nature purement arithmétique, le Lemme \ref{adamslem} 
(déjà remarqué par Adams dans \cite{Adams-J(X)}), montre alors que l'image de 
$\rc_t(\H^j_\dR(X/Y))$ dans ce groupe est annulée par $\N_t$. 
Dans le cas de (a) et (b), cette invariance est démontrée par 
réduction modulo un nombre premier et est une conséquence 
de l'existence de la suite spectrale conjuguée logarithmique en 
caractéristique positive. Le relèvement à la caractéristique nulle 
est rendu possible par un théorème général de changement de base 
en cohomologie étale (Théorème \ref{base-change}) pour 
(a), par l'existence de la variété d'Albanese pour (b) et par 
l'invariance par spécialisation de l'application $\lambda^t$ pour (g). 
Dans le cas de (d), l'invariance est une conséquence 
du théorème d'Adams-Riemann-Roch (cf. la section \ref{adams-sec}) 
et dans le cas de (f), l'invariance provient d'une interprétation 
galoisienne des opérations d'Adams sur l'anneau des représentations 
rationnelles d'un groupe fini. Le point (c) est une conséquence du point 
(b) et du théorème de Rojtman (voir plus bas pour les détails) et le point (e) 
est une conséquence du point (d). Le point (h) est une conséquence du 
point (g) et d'un théorème de Colliot-Thélène, Sansuc et Soulé 
(voir plus bas). 

{\bf Remarque. (1)}. Fixons un corps 
de base de caractéristique nulle $L$. Dans l'article \cite{Ekedahl-VDG}, on considère 
le premier fibré de Gauss-Manin $\H^1_\dR(\CA_L/A_{g,L})$ sur $A_{g,L}$, 
où $A_g$ est le champ classifiant les variétés abéliennes principalement polarisées.  
On y introduit la notion de groupe 
de Chow d'un champ ainsi que la notion de $t$-ème classe de Chern 
à valeurs dans ce groupe. Une conséquence de 
\cite[Prop. 4.2 et Rem. 4.3]{Ekedahl-VDG} est que 
l'ordre de la $t$-ème classe de Chern du fibré $\H^1_\dR(\CA_L/A_{g,L})$ 
dans le groupe de Chow de $A_{g,L}$ est divisible par 
$N_t/2$, si $t$ est pair et si $t\leqslant g$. En particulier, 
si $g=2$,  l'élement $\rc_2(H^1_\dR(\CA_L/A_{g,L}))$ est d'ordre fini 
divisible par $12=N_2/2$. Par ailleurs 
l'équation \refeq{j1} énonce que $8\cdot c_2(H^1_\dR(X/Y))=0$ 
pour toute fibration en surfaces abéliennes sur un base 
lisse et quasi-projective sur $L$. Cela suggère que 
le quotient entre $8$ et $12$ représente une obstruction 
à descendre au champ $A_{g,L}$ l'annulation $8\cdot c_2(H^1_\dR(X/Y))=0$ , qui 
est valable pour toute famille de variétés abéliennes sur un schéma. 

{\bf Remarque. (2)} Les calculs faits dans cet article suggère 
la question suivante, que les auteurs n'osent pas 
élever  au rang de conjecture. Soit $S$ un schéma et soit 
$w:W\to Z$ un morphisme  
de $S$-schémas. On suppose  $W$ et $Z$ munis 
de diviseurs à croisements normaux relativement à $S$ et l'on suppose 
que $w$ est semi-stable relativement à ces diviseurs. 
Par ailleurs, on suppose que la suite spectrale de Hodge 
vers de Rham logarithmique de $W$ sur $Z$ dégénère et que les fibrés  de Gauss-Manin 
logarithmiques $\H^j(W/Z)(\log)$ sont localement libres. 
Il est alors légitime de demander si
\begin{equation}
\psi^p(\H^j(W/Z)(\log))=\H^j(W/Z)(\log)
\label{conj-egm}
\end{equation}
dans $\K_0(Z)[{1\over p}]$, pour presque tout nombre premier $p$. 
Une réponse positive à cette question entrainerait la 
Conjecture \ref{mainconj} (voir la démonstration du Théorème 
\ref{mainth}).  On peut aussi spéculer que 
l'équation \refeq{conj-egm} et donc la Conjecture 
\ref{mainconj} est vérifiée si l'on remplace 
$W/Z$ par un \og motif relatif à singularités 
semi-stable sur $Z$\fg, chaque fois que l'on peut donner un 
sens à l'expression entre guillemets. 

La structure de l'article est la suivante. La section \ref{rappels-sec} 
est consacrée à des rappels sur les morphismes semi-stables, 
les résultats d'Illusie et de Gabber sur la suite spectrale 
de Hodge vers de Rham logarithmique, le théorème 
d'Adams-Riemann-Roch et les classes de Chern en cohomologie étale. 
Dans la section \ref{dem-sec}, on démontre les points du 
Théorème \ref{mainth} dans l'ordre alphabétique. 

{\bf Notations.} Si $G$ est un groupe, on notera 
$\Tor(G)$ le sous-ensemble de $G$ constitué des éléments 
d'ordre fini de $G$. Si $T$ est un espace topologique noethérien, 
on notera $\Irr(T)$ l'ensemble de ses composantes irréductibles. 
Si $H$ est un ensemble et $\CP$ une propriété, l'expression 
\og $\CP(h)$ est vérifiée pour presque tout $h\in H$\fg\ signifie 
que $\CP(h)$ est vérifiée pour tout élément $h\in H'\subseteq H$, 
où $H'$ est un sous-ensemble de $H$ tel que $H\backslash H'$ est fini. 
Si $T$ est l'espace sous-jacent à un schéma et $S\subseteq T$ est un 
sous-ensemble, on écrira $\Zar(S)$ ou $\Zar_T(S)$ pour 
la fermeture de $S$ pour la topologie de Zariski. 

\smallskip
{\bf Remerciements.} 
 Nous remercions P. Brosnan, J.-L. Colliot-Thélène, P. Colmez,  F. Han et J.-P. Serre pour 
 des conversations intéressantes et pour des indications bibliographiques. 
 Nous sommes particulièrement redevables à D. Zagier, pour ses 
 explications sur le Lemme \ref{adamslem} et son contexte arithmétique, à 
 H. Esnault pour ses encouragements et ses explications sur la 
 variété d'Albanese et à R. Pink pour nous avoir encouragé à rédiger la 
 remarque (2) et pour nombre de conversations intéressantes. 
  
 \section{Rappels}

\label{rappels-sec}
 
\subsection{Morphismes semi-stables et suite spectrale conjuguée à pôles logarithmiques}

\label{morss}

Soit $S$ un schéma. Soit $Z$ un schéma lisse sur $S$. 
Soit $D\hookrightarrow Z$ un sous-schéma fermé. 
On dit que $D$ est {\it un diviseur à croisements normaux} dans 
$Z$ relativement à $S$ (cf. \cite[section 1]{Illusie-Reduction}) si pour tout point $z\in Z$, il existe un ouvert $U$ de $Z$ contenant $z$, des
nombres $m,k\in\mN$ avec $k\leqslant m$ et un $S$-morphisme étale $r:U\to\mA^m_S$ tels que 
l'idéal de $D\cap U$ est l'image 
réciproque par $r$ de l'idéal $x_1x_2\cdots x_k$. 
Ici 
$x_1,\dots,x_k$ sont les coordonnées naturelles sur $\mA^n_S$. 
On dira que le quadruplet  $(U,r,m,k)$ est adapté à $D$. 

On notera que dans la définition ci-dessus, $k$ peut être nul. Dans ce 
cas-là, $D\cap U$ est vide. 

Soit $j:X^*:=X\backslash D\hookrightarrow X$ l'inclusion naturelle. 
Le faisceau $\Omega^1_{Z/S}(\log\ D)$ 
est un sous-faisceau en $\CO_Z$-modules de $j_*\Omega^1_{X^*/S}$. 
Il est déterminé  
de manière unique par les conditions suivantes. 

Soit $Z'$ un autre schéma lisse sur $S$ et $D'\hookrightarrow Z'$ un 
diviseur à croisements normaux relativement à $S$. Soit $j':Z'\backslash D'\hookrightarrow Z'$ l'inclusion naturelle. 
\begin{itemize}
\item[$\bullet$]  Si $l:Z\to Z'$ est un $S$-morphisme étale tel 
que $l^*D'=D$, alors il existe une unique 
flèche $\dashrightarrow$ telle que le diagramme
$$
\xymatrix{
l^*\Omega^1_{Z'/S}(\log\ D')\ar[r]\ar@{-->}[d] & l^*j'_*\Omega^1_{Z'/S}\ar[d]^{\sim}\\
\Omega^1_{Z/S}(\log\ D)\ar[r]& j_*\Omega^1_{Z/S}
}
$$
commute et que la flèche $\dashrightarrow$ soit un isomorphisme. 
\item[$\bullet$] Si $l:Z'\to Z$ est un $S$-morphisme étale tel 
que $l^*D=D'$, alors il existe une unique 
flèche $\dashrightarrow$ telle que le diagramme
$$
\xymatrix{
l^*\Omega^1_{Z/S}(\log\ D)\ar[r]\ar@{-->}[d] & l^*j_*\Omega^1_{Z/S}\ar[d]^{\sim}\\
\Omega^1_{Z'/S}(\log\ D')\ar[r]& j'_*\Omega^1_{Z'/S}
}
$$
commute et que la flèche $\dashrightarrow$ soit un isomorphisme. 
\item[$\bullet$] Si $Z=\mA^m_S$ et $D=x_1x_2\cdots x_k$, alors 
$\Omega^1_{Z/S}(\log\ D)$ est libre sur $\mA^m_S$, de base 
$\d x_1/x_1,\dots,\d x_k/x_k, x_{k+1},\dots,x_m$. 
\end{itemize}
En particulier, le faisceau $\Omega^\bullet_{Z/S}(\log\ D)$ est localement libre. 

Supposons maintenant que les composantes irréductibles $D_i$ de $D$ sont toutes 
lisses. On dispose alors d'une suite exacte
$$
0\to\Omega^1_{Z/S}\to\Omega^1_{Z/S}(\log\ D)\stackrel{\Res}{\to}\oplus_i\CO_{D_i}\to 0.
$$
Le morphisme Res s'appelle 
{\it résidu de Poincaré}. Si $Z=\mA^m_S$ et $D=x_1x_2\cdots x_k$, alors 
$$
\Res\big(\sum_{i=1}^k\alpha_i{\d x_i\over x_i}+\sum_{j=k+1}^m\alpha_j x_j\big)=
\oplus_i\alpha_i(\mod\ x_i).
$$

Soit $W$ un autre schéma lisse sur $S$ et $E\hookrightarrow W$ 
un diviseur à croisements normaux relativement à $S$. Soit $g:Z\to W$ un $S$-morphisme. 
On dit que $g$ est {\it semi-stable} relativement à $D$ et $E$ (cf. 
\cite[section 1]{Illusie-Reduction}), s'il 
y a une égalité ensembliste $D=g^{-1}(E)$ et si pour tout point $z\in Z$, 
il existe des nombres $m,n,k,p\in\mN$ et un diagramme
commutatif de $S$-schémas

\xymatrix{
 & & & & & & U\ar[r]^{r}\ar[d]^{g_U} & \mA^m_S\ar[d]^\sigma\\
& & & & & & V\ar[r]^l & \mA^n_S
}

tels que
\begin{itemize}
\item[$\bullet$] $z\in U$;
\item[$\bullet$] le quadruplet $(U,r,m,k)$ (resp. $(V,l,n,p)$) est adapté à 
$D$ (resp. $E$); 
\item[$\bullet$] le morphisme $g_U$ est la restriction de $g$ à $U$ et 
$g(U)\subseteq V$;
\item[$\bullet$] le morphisme $\sigma$ a la forme 
\begin{eqnarray*}
\lefteqn{\sigma(x_1,\dots,x_m)=}\\
&&(x_1x_2\cdots x_{m_1},x_{m_1+1}x_{m_1+2}\cdots x_{m_2},
\dots,x_{m_{p-1}+1}\cdots x_{m_p},x_{m_p+t_1},x_{m_p+t_2},x_{m_p+t_3}\dots x_{m_p+t_{n-p}})
\end{eqnarray*}
où $m_p=k$ et $1\leqslant t_1\leqslant t_2\leqslant t_3\leqslant\dots\leqslant t_{n-p}\leqslant m-k$. 
\end{itemize}

On remarquera qu'une conséquence de cette définition est l'égalité schématique $g^*E=D$. 
On en déduit aussi que $g$ est plat et localement d'intersection complète 
(cf. \cite[Par. 1.3, p. 144]{Illusie-Reduction}). Par ailleurs,  
on voit que les fibres de $g$ sont géométriquement réduites. 

Une autre conséquence est qu'il existe une unique flèche $\dashrightarrow$ telle que 
le diagramme
$$
\xymatrix{
g^*\Omega^1_{W/S}\ar[d]^{g^*}\ar[r]\ar[d] & g^*\Omega^1_{W/S}(\log\ E)\ar@{-->}[d]\\
\Omega^1_{Z/S}\ar[r] & \Omega^1_{Z/S}(\log\ D)
}
\label{useful-diag}
$$
commute. 

Soit $j:U\to Z$ l'ouvert de lissité de $g$. 
Le {\it complexe de de Rham  de $Z$ sur $W$ à pôles logarithmiques 
le long de $D$}, noté $\Omega^\bullet_{Z/W}(\log\ D/E)=\Omega^\bullet_{Z/W}(\log)=\omega^\bullet_{Z/W}$ 
est l'image $j_*\Omega^\bullet_{U/W}$ par $j$ du complexe 
de de Rham de $U$ sur $Z$. 

Au sujet de ce complexe, on peut démontrer les faits suivants. Tout d'abord, 
on dispose d'une suite exacte
$$
0\to g^*\Omega^1_{W/S}(\log\ E){\to} \Omega^1_{Z/S}(\log\ D)\to\omega^1_{Z/W}\to 0
$$
où la deuxième flèche correspondant à la flèche en traitillé $-->$ dans 
le diagramme \refeq{useful-diag}. Enfin, il existe un isomorphisme canonique 
$\Lambda^i(\omega^1_{Z/W})\simeq\omega^i_{Z/W}$.

Soit maintenant $p$ un nombre premier. Supposons que 
$S$ est le spectre d'un corps de caractéristique $p$. 
Considérons le diagramme
$$
\xymatrix{
Z \ar[r]^{F_{Z/W}}\ar[dr]^{g}\ar@/^3pc/[rr]^{F_Z} & Z'\ar[r]^{J}\ar[d]^{g'} & Z\ar[d]^g\\
                            & W\ar[r]^{F_W} & W
}
$$

où $F_W$ (resp. $F_Z$) est le morphisme de Frobenius absolu de $W$ (resp. $Z$). 
Le carré de ce diagramme est par hypothèse cartésien et le morphisme 
$F:=F_{Z/W}$ (appelé morphisme de Frobenius relatif) est l'unique morphisme 
rendant le diagramme commutatif. 

\begin{theor}[Gabber, Illusie]
Les différentielles du  complexe $F_*\omega^\bullet_{Z/W}$ sont 
$\CO_{Z'}$-linéaires et pour tout $i\geqslant 0$, il existe un isomorphisme canonique 
de faisceaux en $\CO_{Z'}$-modules 
$$
C^{-1}:J^*\omega^i_{Z/W}\simeq\CH^i F_*\omega^\bullet_{Z/W}.
$$
\label{th-gabber}
\end{theor}
Remarquons maintenant que, dans la catégorie $\D(W)$ dérivée de 
celle des faisceaux en $\CO_W$-modules, on a 
\begin{eqnarray*}
\R g_*(\omega^\bullet_{Z/W})=\R g'_*(\R F_*(\omega^\bullet_{Z/W})). 
\end{eqnarray*}
Puisque $F$ est fini, 
on dispose d'un isomorphisme $\R F_*(\omega^\bullet_{Z/W})\simeq F_*\omega^\bullet_{Z/W}$ 
dans $\D(Z')$. La deuxième suite spectrale d'hypercohomologie appliquée 
à $F_*\omega^\bullet_{Z/W}$ s'exprime
$$
E_2^{pq}=\R^p g'_*\CH^q(F_*\omega^\bullet_{Z/W})
\Longrightarrow \R^{p+q}g_*(\omega^\bullet_{Z/W})
$$
ou encore, via l'isomorphisme du Théorème \ref{th-gabber}, 
$$
E_2^{pq}=\R^p g'_*J^*\omega^q_{Z/W}
\Longrightarrow \R^{p+q}g_*(\omega^\bullet_{Z/W}).
$$
Comme $F_W$ est plat, cette dernière suite spectrale donne lieu 
à la suite spectrale 
$$
E_2^{pq}=F_W^*\R^p g_*\omega^q_{Z/W}
\Longrightarrow \R^{p+q}g_*(\omega^\bullet_{Z/W}).
$$
On appelle cette dernière suite spectrale {\it la suite 
spectrale conjuguée logarithmique}, ou simplement 
{\it la suite spectrale conjuguée} de $g$. 

On suppose maintenant à nouveau que $S$ est un schéma de base général. 
La première suite spectrale d'hypercohomologie appliquée au 
complexe de de Rham à pôles logarithmiques le long de $D$ 
s'exprime 
$$
E_1^{pq}=\R^q g_*\omega_{Z/W}^p\Longrightarrow \R^{p+q}g_*(\omega^\bullet_{Z/W}).
$$
Cette dernière suite se spécialise en la suite spectrale de Hodge vers de Rham habituelle, lorsque 
$g$ est lisse. Nous l'appellerons {\it suite spectrale de Hodge vers de Rham logarithmique}, ou simplement {\it suite spectrale de Hodge vers de Rham }.  

Une simple comparaison de rangs montre que si la  
suite spectrale de Hodge vers de Rham dégénère en $E_1$, alors  la suite 
spectrale conjuguée (lorsqu'elle existe) dégénère en $E_2$.

\begin{theor}[Illusie]
Si $S$ est le spectre d'un corps de charactéristique nulle, alors la suite spectrale de Hodge vers de Rham dégénère en $E_1$ et les faisceaux 
$\R^q g_*\omega^p_{Z/W}$ sont localement libres.
\label{ill-deg}
\end{theor}

\subsection{Classes de Chern étales}

Dans cette sous-section, on rappelle la définition et les propriétés principales 
des classes 
de Chern en cohomologie étale.
Soit $Z$ un schéma et $n$ un nombre entier inversible 
sur $Z$.  On dispose d'un complexe de schémas en groupes commutatifs 
$$
0\ra \mu_n\ra\mG_m\ra\mG_m\ra 0
$$
qui donne lieu à une suite exacte de faisceaux abéliens 
sur le petit site étale de $Z$, appelée 
 la suite de Kummer. La suite exacte longue de 
cohomologie étale de cette suite contient en particulier 
le morphisme 
$$
\H^1_\et(Z,\mG_m)\ra \H^2_\et(Z,\mu_n).
$$
Par ailleurs, il existe un isomorphisme canonique 
 $\Pic(Z)\simeq H^1_\et(Z,\mG_m)$. On obtient donc un morphisme 
 $$
 \rc^\et_1:\Pic(Z)\ra  H^2_\et(Z,\mu_n).
 $$
 Voir \cite{Grothendieck-Classes} pour tout ceci.

\begin{propdef}
Il existe une unique famille d'opérations 
$\rc^\et_i$ ($i\in\mN$) avec les propriétés suivantes:
\begin{itemize}
\item[\rm\bf (a)] l'opération $\rc^\et_i$ associe un élément $c_i(E)$ 
de $\H^{2i}_\et(Z,\mu_n^{\otimes i})$ à tout faisceau cohérent 
localement libre sur un schéma $Z$ où $n$ est inversible;
\item[\rm\bf (b)] si $g:Z'\ra Z$ est un morphisme de schémas 
 alors 
$f^*\rc^\et_i(E)=\rc^\et_i(f^*E)$ ($n$ étant supposé inversible sur $Z$ et $Z'$);
\item[\rm\bf (c)] l'opération $\rc^\et_1$ est définie comme plus haut;
\item[\rm\bf (d)] si 
$$
0\ra E'\ra E\ra E''\ra 0
$$
est une suite exacte de faisceaux cohérents localement libres sur 
un schéma $Z$ alors on a la relation
$$
\rc^\et_i(E)=\sum_{j+k=i}\rc^\et_j(E')\rc^\et_k(E'')
$$
dans l'anneau de cohomologie $\bigoplus_{m\geqslant 0 } H^{2m}_\et(Z,\mu_l^{\otimes m})$.
\end{itemize}
\label{chernet}
\end{propdef}
Pour la démonstration, voir \cite[Exp. VII, Prop. 3.4]{SGA5}. 

\subsection{Théorèmes de changement de base et de comparaison en cohomologie étale}

Nous rappelons ici deux théorèmes fondamentaux de la théorie de la   
cohomologie étale qui jouerons un rôle essentiel dans la démonstration 
du Théorème \ref{mainth} (a). Soit $n\geqslant 1$. 

\begin{theor}[Deligne et al.]
Soit $g:Z\ra W$ un morphisme de  type fini.  On suppose que 
$n$ est inversible sur $W$ (et donc $Z$). Soit $\cal F$ un faisceau 
 étale constructible en $\mZ/n\mZ$-modules sur $Z$. Alors 
il existe un ouvert dense $U\subseteq W$ tel que 
\begin{description}
\item[(a)] 
pour tout $i\geqslant 0$, $R^i f_*{\cal F}|_U$ est un faisceau 
constructible  sur $U$;
\item[(b)]  pour tout $i\geqslant 0$, le faisceau $R^i f_*{\cal F}|_U$ est invariant par 
tout changement de base au-dessus de $U$.
\end{description}
\label{base-change}
\end{theor}

Pour la démonstration, voir \cite[Finitude, Th. 1.9]{SGA4/2}.  

\begin{theor}[Grothendieck et al.]
Soit $Z$ un schéma séparé et de type fini sur $\mC$. Il existe pour tout $i\geqslant 0$ et tout $n\geqslant 1$ un isomorphisme canonique de comparaison
$$
\kappa_i:H^i_\et(Z,\mZ/n\mZ)\simeq H^i(Z(\mC),\mZ/n\mZ).
$$
\label{comp-iso}
\end{theor}
Pour la démonstration, voir \cite[Exp. XI]{SGA4.3}.

On rappelle que sur $\mC$, le choix d'un isomorphisme de schémas en groupes 
$\mZ/n\mZ\simeq\mu_n$ est équivalent au choix d'une racine $n$-ième primitive 
de l'unité dans $\mC$. 

\begin{theor}[Grothendieck]
Soit $Z$ un schéma lisse sur $\mC$. Identifions 
$\mZ/n\mZ$ et $\mu_n$ via la racine de l'unité 
$\exp(2i\pi/n)$. Alors $\kappa_{2i}(\rc^\et_i(E))$ coïncide 
avec la $i$-ème classe de Chern du fibré vectoriel $E(\mC)$ 
considéré comme fibré vectoriel topologique sur 
$Z(\mC)$. 
\label{compchern}
\end{theor}
Pour la démonstration, voir \cite[Par. 1, p. 243]{Grothendieck-Classes}. 

\subsection{Le théorème d'Adams-Riemann-Roch}

\label{adams-sec}

Soit $W$ un schéma noethérien.
 Nous écrirons $\K_0(W)$ pour le 
groupe de Grothendieck des faisceaux localement libres sur $W$.
Le produit tensoriel induit une application bilinéaire de $\K_0(W)$ dans 
$\K_0(W)$ qui en fait un anneau commutatif unifère. Pour tout 
morphisme de schémas noethériens $g:Z\to W$, le foncteur $f^*$
(image réciproque des faisceaux en $\CO_W$-modules) induit 
un morphisme d'anneaux $g^*:\K_0(W)\to\K_0(W)$. 

On rappelle que pour tout $k\geqslant 0$, l'opération 
d'Adams $\psi^k:\K_0(W)\to\K_0(W)$ est l'unique endomorphisme 
d'anneau compatible aux images réciproques tel que 
$\psi^k(L)=L^{\otimes k}$ pour tout fibré en droites $L$. 

Les opérations d'Adams satisfont les deux relations 
de compatibiités suivantes.

Supposons le temps de la prochaine phrase que $W$ est lisse 
sur un corps. Il existe alors pour tout $t\geqslant 0$ des applications classe 
de Chern $\rc_t:\K_0(W)\to{\rm CH}^t(W)$ et l'on a 
\begin{equation}
\rc_t(\psi^k(w))=k^t\cdot c_t(w).
\label{comp-ad-ch}
\end{equation}
Supposons à nouveau juste le temps de la prochaine phrase
que $p$ est un nombre premier et que $W$ est un schéma sur $\mF_p$ 
tel que le morphisme de Frobenius absolu $F_W:W\to W$ est plat. 
On a alors 
\begin{equation}
F^*_W=\psi^p.
\label{comp-ad-fro}
\end{equation}
Les deux compatibilités sont des conséquences du principe 
de scindage (cf. \cite[Th. 2.1]{Quillen-Higher}). 

Nous aurons aussi besoin des {\it classes cannibales} $\theta^k$. 
Pour chaque $k\in\mN^*$, $\theta^k$ est une opération 
associant des éléments de 
$\K_0(W)$ à des faisceaux localement libres sur $W$ et 
jouissant des propriétés suivantes. Tout d'abord, pour 
tout fibré en droites $L$, on a 
$$
\theta^k(L)=1+L+L^{\otimes 2}+\dots+L^{\otimes (k-1)}.
$$
Par ailleurs, pour toute suite exacte
$$
0\to E'\to E\to E''\to 0
$$
de faisceaux localement libres, on a $\theta^k(E')\otimes\theta^k(E'')=\theta^k(E)$. 
Enfin les classes $\theta^k$ sont compatibles aux images 
réciproques. 
Ces trois propriétés déterminent l'opération $\theta^k$. 

Dans le contexte des classes cannibales, nous mentionnerons le lemme suivant:  
\begin{lemme}
Si $W$ est muni d'un faisceau ample, alors l'élement 
$\theta^k(E)$ est inversible dans 
$\K_0(W)[{1\over k}]$ pour tout faisceau localement libre $E$ sur $W$. 
\end{lemme}
Pour la démonstration voir par exemple \cite[sec. 4, Prop. 4.2]{Rossler-Adams} (faute d'une référence 
canonique). 

On suppose maintenant que $W$ est régulier et possède un fibré 
en droites ample. Soit $\K'_0(W)$ le groupe de Grothendieck des faisceaux 
cohérents sur $Y$. On peut démontrer que l'homomorphisme 
de groupes naturel $\K_0(W)\to\K'_0(W)$ est un isomorphisme. 
Soit maintenant $g:Z\to W$ un morphisme projectif et localement 
d'intersection complète de schémas. On définit alors 
le morphisme de groupes $\R g_*:\K_0(Z)\to\K'_0(W)$ 
par la formule 
$$
\R g_*(E):=\sum_{l\geqslant 0}(-1)^l\R^l g_*(E)
$$
et on notera aussi $\R g_*$ le morphisme 
$\K_0(Z)\to\K_0(W)$ obtenu par composition. 

Par hypothèse, on dispose d'une factorisation 
$$Z\stackrel{i}{\hookrightarrow} P\stackrel{p}{\to} W$$ de $g$ telle que 
$i$ est une immersion régulière et $p$ est lisse. On définit
$$
\theta^k(g)^{-1}:=\theta^k(\Omega_p)^{-1}\otimes\theta^k(N)
$$
où $N$ est le fibré conormal de l'immersion $i$. On 
peut montrer que $\theta^k(g)^{-1}$ ne dépend pas 
de la factorisation.

\begin{theor}[th. d'Adams-Riemann-Roch]
 Pour tout $z\in\K_0(Z)[{1\over k}]$, l'égalité 
$$
\psi^k(\R g_*(z))=\R g_*(\theta^k(g)^{-1}\otimes\psi^k(z))
$$
est vérifiée dans $\K_0(W)[{1\over k}]$.
\label{arr}
\end{theor}
Pour une démonstration complète (reposant de fa\c{c}on essentielle sur 
les idées de \cite{SGA6} et \cite{BFQ}), voir \cite{Koeck-Grothendieck}. 

\section{Démonstration du Théorème \ref{mainth}}

\label{dem-sec}

\subsection{Démonstration de (a)} 

On peut tout d'abord supposer sans restreindre la généralité que $Y$ est connexe.

\begin{lemme} 
Supposons que pour tout nombre premier $l$ et pour tout $n\geqslant 1$, 
l'image de 
 $
 \N_t\cdot \rc_t(\H_\dR^{j}(X/Y))
 $
 dans $\H^{2t}(Y(\mC),\mZ/l^n\mZ))$ est nulle. Alors 
 le point (a) du Théorème \ref{mainth} est vérifié.
 \end{lemme}
\beginProof
Considérons les application naturelles
\begin{center}
\xymatrix{
&&&&\H^{2t}(Y(\mC),\mZ_{l})\ar[r]^{\lambda\,\,\,\,\,\,\,\,\,\,\,\,\,} &
\prolim_{j}\H^{2t}(Y(\mC),\mZ/l^j\mZ)\\
&&&&\H^{2t}(Y(\mC),\mZ)\ar[u]^\mu & 
}
\end{center}
L'espace $X(\mC)$ muni de la topologie ordinaire est de type fini. Cela implique que 
l'application $\lambda$ est un isomorphisme 
(voir \cite[3.3]{Grothendieck-Classes} pour les détails). 
Par ailleurs, compte tenu du fait que $\mZ_l$ est plat sur $\mZ$, cela 
implique aussi que  
$\H^{2t}(Y(\mC),\mZ_{l})\simeq\H^{2t}(Y(\mC),\mZ)\otimes\mZ_l$ 
(au vu de la suite universelle des coefficients pour la cohomologie). 
En utilisant le fait que $\H^{2t}(Y(\mC),\mZ)$ est de type fini, 
on en déduit que le noyau $\ker(\mu)$ de 
$\mu$ est de torsion et que $\ker(\mu)$ ne contient aucun élément non nul 
qui soit de  $l^\infty$-torsion. Autrement dit, le morphisme 
$\ker(\mu)\to\ker(\mu)$ de multiplication par $l$ est injectif.

Remarquons maintenant que $\mu$ et $\lambda$ sont par construction compatibles 
à la formation des classes des Chern. On en déduit que l'image de
 $
 e:=\N_t\cdot \rc_t(\H_\dR^{j}(X/Y)(\log))
 $
 dans $\H^{2t}(Y(\mC),\mZ)$ est de torsion. Par ailleurs, si 
 $e\not=0$, alors pour tout $k\geqslant 1$ et tout 
 nombre premier $l$, $l^k\cdot e\not=0$. Ceci 
 est la contradiction que achève la preuve du lemme. 
 \endProof

 Au vu du dernier lemme et du Théorème \ref{compchern}, il suffit 
 de démontrer que pour tout nombre premier $l$ et pour 
tout $n\geqslant 1$ 
la $t$-ème classe de Chern de $\H^{j}_\dR(X/Y))(\log)$ dans 
$\H^{2t}(Y(\mC),\mu_{l^n})$ est annulée par $N_t$. 

Fixons donc un nombre premier $l$ et un nombre entier 
$n\geqslant 1$. Par abus de notation, jusqu'à la fin de la démonstration de (a), 
nous écrirons 
$c_t(\H^{j}_\dR(X/Y)(\log))$ pour la classe de Chern de 
$\H^{j}_\dR(X/Y)(\log)$ dans $\H^{2t}(Y,\mu_{l^n})\stackrel{\kappa_{2t}}{\simeq}\H^{2t}(Y(\mC),\mZ/l^n\mZ)$.

Soit $L$ un corps de type fini sur $\mQ$ (comme corps) tel que 
le triplet $f:X\to Y$ possède un modèle $f_0:X_0\to Y_0$ sur $L$.  Quitte 
à remplacer $L$ par une de ses extensions finies, on peut supposer 
que $L$ contient les racines $l^n$-èmes de $1$. Soit 
$T$ un schéma affine, intègre, lisse et de type fini sur  
$\mZ$, dont le corps de fonctions est isomorphe à $L$. Quitte 
à remplacer $T$ par l'un de ses ouverts, on peut supposer qu'il existe 
des modèles lisses 
$\widetilde{X}$ et $\widetilde{Y}$ de $X_0$ et $Y_0$ sur $T$ et 
qu'il existe un $T$-morphisme lisse et projectif
$\widetilde{f}:\widetilde{X}\ra\widetilde{Y}$ qui est 
un modèle de $f_0$. Toujours quitte à réduire la taille de $T$, on peut supposer que la suite spectrale de Hodge vers de Rham 
de $\widetilde{f}:\widetilde{X}\ra\widetilde{Y}$ dégénère et que les faisceaux de cohomologie 
de de Rham relative correspondants sont localement libres. Enfin, on peut 
supposer que $l$ est inversible sur $T$. Soit $\wt{p}:\wt{Y}\to T$ 
le morphisme structural. 
Au vu du Théorème \ref{base-change}, on aussi peut encore supposer 
que les faisceaux 
étales $\R^i\widetilde{p}_*\mu_{l^n}$ ($i\geqslant 0$) sont
invariants par changement de base au-dessus de $T$. 
Abbrévions 
$\A^{t}(\bullet):=\H^{2t}_\et(\bullet,\mu_{l^n}^{\otimes 2t})$. 
Soit $\mathfrak p$ un point fermé de $T$. Soit $p$ la charactéristique 
du corps résiduel en $\mathfrak p$ (qui est nécessairement positive). 
Le corps résiduel de $\mathfrak p$ est alors isomorphe 
à un corps fini $\mF_q$, où $q$ est une puissance de $p$. 
Soit $R_{\mathfrak p}$ la 
Henselisation stricte de l'anneau local $\CO_{T,\mathfrak p}$ de 
$T$ en $\mathfrak p$. On dispose par construction d'un diagramme 
commutatif de morphismes
$$
\xymatrix{
\Spec\ \mC\ar[r] & \Spec\ L^\sep\ar[r]\ar[dr] & \Spec\ \CO_{T,{\mathfrak p}}
\ar[r]&T&
& \Spec\ \mF_q\ar[ll]\\
& & \Spec\ R_{\mathfrak p}\ar[u] & \Spec\ \mtr{\mF}_q\ar[l]\ar[urr] &\\
}
$$
(où $L^\sep$ est la clôture séparable de $L$). Par changement de base, 
ce diagramme induit le diagramme commutatif
$$
\xymatrix{
\widetilde{Y}_\mC=Y\ar[r]^{s_1} & \widetilde{Y}_{L^\sep}\ar[r]^{s_2}\ar[dr]^{s_6} & 
\widetilde{Y}_{\CO_{T,{\mathfrak p}}}
\ar[r]^{s_3}&\widetilde{Y}&
& \widetilde{Y}_{\mF_q}\ar[ll]_{s_4}\\
& & \widetilde{Y}_{R}\ar[u]^{s_7} & \widetilde{Y}_{\mtr{\mF}_q}\ar[l]^{s_8}\ar[urr]^{s_5} &\\
}
$$
 Soit 
$\pi:\widetilde{Y}_{\mtr{\mF}_q}\ra\Spec\ \mtr{\mF}_q$ le morphisme de structure. 
On dispose alors de deux suites spectrales: la suite 
spectrale de Hodge vers de Rham logarithmique 
$$
E_1^{rs}=\R^s f_*\Omega^r_{\widetilde{X}_{\mtr{\mF}_q}/\widetilde{Y}_{\mtr{\mF}_q}}
\Longrightarrow H^{r+s}_\dR(\widetilde{X}_{\mtr{\mF}_q}/\widetilde{Y}_{\mtr{\mF}_q})
$$
et la suite spectrale conjuguée logarithmique (voir sous-section \ref{morss}) 
\begin{equation*}
E_2^{rs}=F^*_Y \R^r f_*(\Omega^s_{\widetilde{X}_{\mtr{\mF}_q}/
\widetilde{Y}_{\mtr{\mF}_q}})(\log)\Longrightarrow 
H^{r+s}_\dR(\widetilde{X}_{\mtr{\mF}_q}/\widetilde{Y}_{\mtr{\mF}_q})(\log).
\end{equation*}
Comme la première suite spectrale dégénère par hypothèse, il en est de même 
de la deuxième. En utilisant les compatibilités \refeq{comp-ad-ch} et \ref{comp-ad-fro}, on obtient donc l'égalité 
$$
\rc_t(\H^{j}_\dR(\widetilde{X}_{\mtr{\mF}_q}/\widetilde{Y}_{\mtr{\mF}_q})(\log))=
\rc_t(F^*\H^{j}_\dR(\widetilde{X}_{\mtr{\mF}_q}/\widetilde{Y}_{\mtr{\mF}_q})(\log))=
p^t\cdot \rc_t(\H^{j}_\dR(\widetilde{X}_{\mtr{\mF}_q}/\widetilde{Y}_{\mtr{\mF}_q})(\log))
$$
dans $A^t(\widetilde{Y}_{\mtr{\mF}_q})$, qui est le pivot du présent article.
L'invariance par changement de base assure maintenant que le morphisme 
$$
s^*_8:A^\bullet(\widetilde{Y}_R)\ra A^\bullet(\widetilde{Y}_{\mtr{\mF}_q})
$$
est un isomorphisme.  On déduit que 
$$(1-p^t)\cdot s^*_7 s^*_3(\rc_t(\H^{j}_\dR(\widetilde{X}/\widetilde{Y})(\log)))=0$$ et donc 
que 
\begin{eqnarray}
(1-p^t)\cdot s^*_1 s^*_2 s^*_3(\rc_t(\H^{j}_\dR(\widetilde{X}/\widetilde{Y})(\log)))=
(1-p^t)\cdot \rc_t(\H^{j}_\dR({X}/{Y})(\log))=0
\label{fundeq}
\end{eqnarray}
dans $\A^t(Y)$. 

Remarquons maintenant que l'image de $T$ dans $\Spec\ \mZ$ est constructible par le 
théorème de Chevalley et que chaque fibre de $T$ sur $\Spec\ \mZ$ contient 
un point fermé (car $T$ est de type fini sur $\mZ$). 
L'égalité \refeq{fundeq} est donc vérifiée pour tous les nombres premiers 
$p$ en dehors d'un ensemble fini. Nous allons maintenant 
appliquer le lemme suivant.
\begin{lemme}
Soit $E\subseteq\Spec\ \mZ$ un ensemble de densité de Dirichlet 
$1$. Soit $f:E\ra\mN$ une fonction. Alors le nombre
$$
\pgcd\{(1-p^t)\cdot p^{f(p)}\}_{p\in E}
$$
divise $\N_t$. 
\label{adamslem}
\end{lemme}
\beginProof
Soit $\mu\geqslant 1$. On rappelle les faits suivants. 
Si $p$ est un nombre premier impair, il existe une unité 
d'ordre $(p-1)p^{\mu-1}$ dans $(\mZ/p^\mu\mZ)^*$. Si $p=2$ et 
$\mu>3$ alors il existe une unité d'ordre $p^{\mu-2}$ dans 
$(\mZ/p^\mu\mZ)^*$. Si $\mu=2,3$, il existe une unité d'ordre $2$ 
dans $(\mZ/2^\mu\mZ)^*$. 
Pour une démonstration de ces faits, voir 
\cite[Ch. 4, Prop. 3.32 et Cor. 3.34]{Demazure-Cours}.

Soit maintenant $c:=\pgcd\{(1-p^t)\cdot p^{f(p)}\}_{p\in E}$. 
Pout tout nombre premier $p$, soit $c_p$ (resp. $t_p$) la multiplicité de 
$p$ dans $c$ (resp. $t$). 
Soit $E'$ l'ensemble $E$ privé des nombres premiers 
divisant $c$. Par construction du nombre $c$, pour tout $p\in E'$, 
les $c$-èmes racines de l'unité dans 
$\mtr{\mF}_p$ sont contenues dans le groupe 
cyclique $\mF_{p^t}^*$. En conséquence, 
on a $\mF_p(\mu_c)\subseteq\mF_{p^t}$. On en déduit que 
$[\mF_p(\mu_c):\mF_p]\ |\ t$. 

Ceci implique que 
pour tout $p\in E'$, le symbole d'Artin $\sigma_p\in G:=\Gal(\mQ(\mu_c)|\mQ)$  a la propriété 
suivante:
\begin{equation}
\ord(\sigma_p)\ |\ t
\label{artinid}
\end{equation}
Rappelons que $G\simeq(\mZ/c\mZ)^*\simeq 
\prod_{p|c}(\mZ/p^{c_p}\mZ)^*$. 
Au moyen du théorème de densité de Tchebotareff, on déduit de 
\refeq{artinid} que 
$$
2^{c'_2-2}\prod_{p\not=2,\ p|c}(p-1)p^{c_p-1}|t
$$
où $c'_2=c_2$ si $c_2>3$, $c'_2=4$ si $c_2=2,3$ et 
$c'_2=2$ si $c_2=0,1$. 
En particulier, si $p$ est impair et $p|c$, alors $p-1|t$ et 
$c_p\leqslant t_p+1$. 
De même, on a $c'_2\leqslant t_2+2$. On en déduit que $c\ |\ \N_t$. 
\endProof
Le Lemme \ref{adamslem} est une variation d'un résultat d'Adams 
démontré dans \cite{Adams-J(X)}. 

Soit maintenant $E$ l'image de $T$  dans 
$\Spec\ \mZ$ et $f:E\to\mN$ la fonction constante de valeur $0$. 
Le Lemme \ref{adamslem} conclut la démonstration de (a). 

\subsection{Démonstration de (b) \& (c)}

\label{ssec-bc}

On peut encore supposer, sans restriction de généralité, que 
$Y$ est connexe. 

{\it Commen\c{c}ons par la démonstration de (b)}. 
Soit encore une fois $L$ un corps de type fini sur $\mQ$ (comme corps) tel que 
le triplet $f:X\to Y$ possède un modèle $f_0:X_0\to Y_0$ sur $L$.  
Soit $\kappa_L:=\sum_{i}n_i P_i$ une combinaison linéaire formelle de points fermés de $X_0$ 
dans $X_0({L})$ telle que l'image de $\kappa_L$ dans 
${\rm CH}^{d_Y}(Y_0)$ est égale à $\rc_{d_Y}(\H^j_\dR(X_0/Y_0)(\log))$. 
Quitte à remplacer $L$ par une de ses extensions finies, une telle combinaison linéaire existe. 
Soit 
$T$ un schéma affine, intègre, lisse et de type fini sur  
$\mZ$, dont le corps de fonctions est isomorphe à $L$. Comme avant, quitte 
à remplacer $T$ par l'un de ses ouverts, on peut supposer qu'il existe 
des modèles lisses 
$\widetilde{X}$ et $\widetilde{Y}$ de $X_0$ et $Y_0$ sur $T$ et 
qu'il existe un $T$-morphisme lisse et projectif
$\widetilde{f}:\widetilde{X}\ra\widetilde{Y}$ qui est 
un modèle de $f_0$. 
Toujours quitte à réduire la taille de $T$, on peut supposer que la suite spectrale de Hodge vers de Rham 
de $\widetilde{f}:\widetilde{X}\ra\widetilde{Y}$ dégénère et que les faisceaux de cohomologie 
de de Rham relative correspondants sont localement libres. Enfin, on peut 
supposer que les fibres géométriques de $Y_0$ sur $T$ sont connexes et projectives. 

Soit $\kappa:=\sum_{i}n_i\cdot\Zar(P_i)$ la combinaison linéaire formelle
des clôtures de Zariski des $P_i$ dans $T$. Quitte à restreindre $T$, on 
peut supposer que les $\Zar(P_i)$ sont les images de 
sections $\widetilde{P}_i\in\widetilde{Y}(T)$. 
La construction des classes de Chern montre qu'il existe un ouvert
$U\subseteq T$ tel que l'image de $\sum_{i}n_i\cdot\widetilde{P}_i(u)$ dans 
${\rm CH}^{d_Y}(Y_{0,u})$ est égale à $\rc_{d_Y}(\H^j_\dR(X_{0,u}/Y_{0,u})(\log))$ 
pour tout $u\in U$. On remplace $T$ par $U$. 

On munit maintenant le schéma $Y$ d'un point de base $P\in Y(L)$, dont on peut 
à nouveau supposer qu'il s'étend en une section $\wt{P}\in\wt{Y}(T)$. 
Nous écrirons $\Pic_{\wt{Y}/T}$ pour le schéma sur $\widetilde{Y}$ 
représentant le foncteur sur les $T$-schémas $S$

\medskip
$S\mapsto\{\CL|\CL\ $fibré en droites sur $S\times_T\wt{Y}$ rigidifiés le long de 
$S\times_T\wt{P}\}$
\medskip

Outre le fait qu'il existe, il est montré dans que \cite[Th. 3.1]{FGA-232} que ce schéma est séparé et localement de type fini sur $T$. 
Il est aussi montré au même endroit que $\Pic_{\wt{Y}/T}$ est une réunion 
disjointe de sous-schémas ouverts et fermés, qui sont quasi-projectifs au-dessus de $U$. 
Ainsi, quitte à restreindre $T$ encore une fois,  on peut supposer qu'il existe 
un $T$-sous-schéma en groupes $\Pic^0_{\widetilde{Y}/T}$ ouvert et fermé de 
$\Pic_{\widetilde{Y}/T}$ tel que 
$\Pic^0_{\widetilde{Y}/T,s}=\Pic^0_{\widetilde{Y}_s/T_s}$ 
pour tout $s\in T$ 
et tel que $\Pic^0_{\widetilde{Y}/T}$ est un schéma abélien. 
On rappelle que $\Pic^0_{\widetilde{Y}_s/T_s}$ est la composante 
neutre du schéma $\Pic_{\widetilde{Y}_s/T_s}$ et que 
 $\Pic^0_{\widetilde{Y}_s/T_s}$ est ouvert 
et fermé dans $\Pic_{\widetilde{Y}_s/T_s}$. 
On rappelle aussi que $\Pic^0_{\widetilde{Y}_s/T_s}$ est de type fini et 
géométriquement irréductible sur $\kappa(s)$ 
(cf. \cite[Exp. $\rm VI_A$, Par. 2.3]{SGA3-1}). 

On considère maintenant le morphisme $\wt{\alb}:\wt{Y}\to(\Pic^0_{\widetilde{Y}/T})^\vee$ défini par la restriction du fibré universel $\CP$ sur 
$\Pic_{\widetilde{Y}/T}\times_T\wt{Y}$, muni de 
sa rigidification naturelle $\CP|_{0\times_T\wt{Y}}$.  Par construction, 
$\wt{\alb}_s$ coïncide avec le morphisme 
d'Albanese $\wt{Y}_s\to\Alb(\wt{Y}_s)\simeq(\Pic^0_{\wt{Y}_s/T_s})^\vee$ 
pour tout $s\in T$. 

Soit maintenant $\mathfrak p$ un point de $T$ dont le corps résiduel est
 isomorphe 
à un corps fini $\mF_q$, où $q$ est une puissance de $p$. Soit 
$\pi:\widetilde{Y}_{\mtr{\mF}_q}\ra\Spec\ \mtr{\mF}_q$ le morphisme de structure. 
On dispose alors comme avant de deux suites spectrales: la suite 
spectrale de Hodge vers de Rham logarithmique 
$$
E_1^{rs}=\R^s f_*\Omega^r_{\widetilde{X}_{\mtr{\mF}_q}/\widetilde{Y}_{\mtr{\mF}_q}}
\Longrightarrow H^{r+s}_\dR(\widetilde{X}_{\mtr{\mF}_q}/\widetilde{Y}_{\mtr{\mF}_q})
$$
et la suite spectrale conjuguée logarithmique 
\begin{equation*}
E_2^{rs}=F^*_Y \R^r f_*(\Omega^s_{\widetilde{X}_{\mtr{\mF}_q}/
\widetilde{Y}_{\mtr{\mF}_q}})(\log)\Longrightarrow 
H^{r+s}_\dR(\widetilde{X}_{\mtr{\mF}_q}/\widetilde{Y}_{\mtr{\mF}_q})(\log).
\end{equation*}
Comme la première suite spectrale dégénère par hypothèse, il en est de même 
de la deuxième. On obtient donc l'égalité
$$
(1-p^{\rc_{d_Y}})\cdot \wt{\alb}(\rc_{d_Y}(\H^j_\dR(\wt{X}_\mfp/\wt{Y}_\mfp)(\log)))=0
$$
dans $\Alb(\wt{Y}_\mfp)({\mF}_q)$. Soit maintenant un point fermé 
$Q\in T_\mQ$. Son corps résiduel $\kappa(Q)$ est 
une extension finie de $\mtr{\mQ}$. La clôture de Zariski 
$\Zar(Q)$ de $Q$ dans $T$ contient alors un ouvert 
$V$ qui est isomorphe à un ouvert de 
$\Spec\ \CO_{\kappa(Q)}$, où 
$\CO_{\kappa(Q)}$ est la clôture intégrale de $\mZ$ dans 
$\kappa(Q)$ (utiliser \cite[IV, Cor. 6.12.6]{EGA} et 
le \og Main Theorem\fg\ de Zariski, cf. par ex. \cite[4.4, Cor. 4.6]{Liu-Algebraic}). 
Par ailleurs, tous les corps résiduels de tous les points de 
$V\backslash\{Q\}$ sont des corps finis. On obtient 
ainsi une variété abélienne $A:=\Alb(\wt{Y}_Q)$ sur le corps de 
nombres $\kappa(Q)$ et un point $a\in A(\kappa(Q))$, donné par 
$\alb(\rc_{d_Y}(\H^j_\dR(\wt{X}_Q/\wt{Y}_Q)$, avec la propriété 
suivante. Pour presque toutes les places finies 
$\mfp$ de  $\CO_{\kappa(Q)}$, la réduction modulo $\mfp$ de $a$ 
est dans le noyau de $1-p^{\rc_{d_Y}}$. Ici $p$ est la caractéristique 
résiduelle de $\mfp$. 

Le théorème de Pink \cite[Th. 4.7]{Pink-Order} implique alors 
que $a$ est un point de torsion. 
Par ailleurs, comme $a$ est un point de torsion, il est par construction 
annulé par le nombre $\pgcd\{1-p^{\rc_{d_Y}}\}_{p\in I}$ 
où $I$ est un ensemble cofini dans $\Spec\ \mZ$. Le Lemme 
\ref{adamslem} implique alors que $\N_{d_Y}\cdot a=0$, ou autrement dit, que 
$$
\N_{d_Y}\cdot \wt{\alb}(\rc_{d_Y}(\H^j_\dR(\wt{X}_Q/\wt{Y}_Q)(\log)))=0
$$
dans $\Alb(\wt{Y}_\mfp)(\mtr{\kappa(Q)})$. Comme 
les points fermés sont denses dans $T_\mQ$, on conclut que 
$$
\N_{d_Y}\cdot \wt{\alb}(\rc_{d_Y}(\H^j_\dR(\wt{X}_L/\wt{Y}_L)(\log)))=0
$$
dans $\Alb(\wt{Y}_\mfp)({L})\subseteq \Alb(\wt{Y}_\mfp)(K)$, ce qu'il fallait démontrer.

{\it Passons maintenant à la démonstration de (c)}. 

On rappelle que le théorème de Rojtman \cite{Rojtman-torsion} implique que la 
l'application
$
\alb:{\rm CH}^{d_Y}(Y)\to\Alb(Y)(K)
$
induit un isomorphisme $\Tor({\rm CH}^{d_Y}(Y))\to\Tor(\Alb(Y)(K))$ des sous-groupes 
de torsion correspondants. 
Le point (c) est une conséquence de ce résultat et du point (b). 

\subsection{Démonstration de (d) \& (e)}

{\it Commen\c{c}ons par le point (d).} 

Les calculs faits dans cette sous-section sont inspirés par les 
calculs faits dans \cite{Esnault-Viehweg-Chern}. 
Nous allons utiliser la notation de la sous-section \ref{morss}. 
Appliquons le théorème d'Adams-Riemann-Roch 
au complexe de de Rham de $f$ à pôles logarithmiques 
le long de $D$. On obtient l'égalité
\begin{equation}
\psi^k(\R f_*(\Lambda_{-1}(\Omega^1_{X/Y}(\log)))=
\R f_*(\theta^k(f)^{-1}\otimes\psi^k(\Lambda_{-1}(\Omega^1_{X/Y}(\log)))).
\label{eqadams}
\end{equation}
dans $\K_0(Y)[{1\over k}]$. 
On souhaite maintenant expliciter le terme
\begin{equation}
\theta^k(f)^{-1}\otimes\psi^k(\Lambda_{-1}(\Omega^1_{X/Y}(\log))).
\label{probterm}
\end{equation}
Soit $\{E_u\}$ l'ensemble des composantes irréductibles 
de $E$ et $\{D_i\}$ l'ensemble des composantes irréductibles 
de $D$. Soit $\CO_{E_0}:=\oplus_u\CO_{E_u}$ et 
$\CO_{D_0}:=\oplus_i\CO_{D_i}$. 
Considérons le diagramme commutatif
\begin{equation}
\begin{matrix}
& & 0 & & 0 & & 0\cr
& &\downarrow &\hskip-39pt{\ss (1)} & \downarrow & & \downarrow & \hskip-70pt{\ss (2)}\cr
0 & \stackrel{(3)}{\rightarrow} & f^*\Omega_Y & {\ra} & f^*\Omega_Y(\log\ E) & 
\stackrel{f^*{\rm Res}_Y}{\ra} & f^*{\cal O}_{E_0} & \ra & 0\cr
&  & \downarrow & & \downarrow & & \downarrow&\hskip-70pt{\ss (4)} &\cr
0 & \ra &  \Omega_X & \ra & \Omega_X(\log\ D)& \stackrel{{\rm Res}_X}{\ra} & {\cal O}_{D_0} &  \to 
& 0\cr
& & \downarrow & & \downarrow& &\downarrow& &\cr
0 &\stackrel{(5)}{\to} & \Omega_{X/Y} & \ra & \Omega^1_{X/Y}(\log)& \ra & 
 {\cal O}_{D_0}/f^*{\cal O}_{E_0} &\stackrel{(6)}{\to} & 0\cr
 & & \downarrow & & \downarrow& &\downarrow& &\cr
 &  & 0 & & 0 & & 0\cr
 \end{matrix}
 \label{fundmatr}
\end{equation}
Seules les flèches numérotées requièrent une justification. 
La flèche (1) est justifiée par le fait que le morphisme $f$ 
est génériquement lisse et par le fait que $X$ et $Y$ sont lisses
 sur $K$. La flèche 
(3) est justifiée par le fait que le morphisme $f$ est plat. La flèche 
(4) est l'unique flèche rendant le diagramme commutatif. 
Cette définition de (4) peut être explicitée de la manière suivante. 
 Soit $u$ et $i$ tels que 
$D_i\subseteq f^{-1}E_u$. Du fait que 
$D_i$ est une composante irréductible de $f^{-1}E_u$, 
on dispose d'un morphisme 
$$
a_{u,i}:\CO_{f^*\CO_{E_u}}\to \CO_{D_i}.
$$ 
La flèche (4) est alors  
décrite par la formule
$$
\oplus_u e_u\mapsto\sum_{u}\sum_{i,D_i\subseteq f^{-1}E_u}(0,\dots,a_{u,i}(e_u),\dots,0)
$$
où le terme $a_{u,i}(e_u)$ est d'indice $i$. 
Comme $g^*E_u$ est un schéma réduit pour tout $u$ par définition de la 
semi-stabilité, le morphisme $a_{u,i}$ est injectif pour tout couple 
$u,i$ tel que 
$D_i\subseteq g^{-1}E_u$. Conséquemment, la flèche (4) 
est injective, d'où la flèche (2). 
La flèche (5) est maintenant justifiée par le lemme du serpent. 
La flèche (6) est justifiée par une chasse au diagramme élémentaire.

Revenons maintenant au terme \refeq{probterm}. Un examen 
de la troisième ligne du diagramme \refeq{fundmatr} donne
\begin{eqnarray*}
\theta^k(f)^{-1}&=&
\theta^k(\Omega^1_{X/Y}(\log))^{-1}\otimes\theta^k({\cal O}_{D_0}/f^*{\cal O}_{E_0}).
\end{eqnarray*}
et donc le terme \refeq{probterm} peut s'exprimer comme 
\begin{equation}
\theta^k(\Omega_{X/Y}(\log))^{-1}\otimes\theta^k({\cal O}_{D_0}/f^*{\cal O}_{E_0})
\otimes\psi^k(\Lambda_{-1}(\Omega^1_{X/Y}(\log))).
\label{probtermp}
\end{equation}
Nous allons démontrer l'équation
\begin{equation}
\theta^k({\cal O}_{D_0}/f^*{\cal O}_{E_0})\otimes\psi^k(\Lambda_{-1}(\Omega^1_{X/Y}(\log)))=\psi^k(\Lambda_{-1}(\Omega^1_{X/Y}(\log)))
\label{pivotadams}
\end{equation}
qui est le nerf de notre calcul. 

Pour ce faire, nous allons d'abord démontrer les trois lemmes suivants.
\begin{lemme}
Soit $D_i, D_j$ des composantes irréductibles  de $D$. Alors pour tout 
sous-schéma fermé $Z$ de $D_i\cap D_j$, il existe un morphisme 
surjectif
$$
\Omega^1_{X/Y}(\log)\to\CO_Z
$$
de faisceaux. 
\label{firstlem}
\end{lemme}
\beginProof
On rappelle le fait suivant. Soit $T$ un schéma  
et $\{T_t\}_{t\in\{0,\dots,q\}}$ un ensemble fini de sous-schémas fermés 
tels que $\cup_{t}T_t=T$. On dispose alors d'un complexe
\begin{eqnarray}
& \CO_T\to\oplus_t\CO_{T_t}\stackrel{d}{\to}\oplus_{t_0<t_1}\CO_{T_{t_0}\cap T_{t_1}}
\label{IDW}
\end{eqnarray}
où 
$$
(d(\alpha))_{t_0,t_1}=
\alpha_{t_0}|_{\CO_{T_{t_0}\cap T_{t_1}}}-\alpha_{t_1}|_{\CO_{T_{t_0}\cap T_{t_1}}}.
$$
On déduit de l'existence de ce complexe et de la description de $d$ l'existence d'une surjection
$$
{\cal O}_{D_0}/f^*{\cal O}_{E_0}\to \CO_{D_i\cap D_j}.
$$
Si l'on compose ce morphisme avec la surjection 
$
 \Omega^1_{X/Y}(\log)\ra 
 {\cal O}_{D_0}/f^*{\cal O}_{E_0}$ à gauche et avec la surjection 
 $\CO_{D_i\cap D_j}\to \CO_Z$ à droite, on obtient la surjection promise. 
\endProof

\begin{lemme}
Soit $z:Z\hookrightarrow X$ un sous-schéma fermé de $X$. Soit 
$V$ un faisceau localement libre sur $X$. On suppose 
qu'il existe une suite exacte
$$
V\to \CO_Z\to 0.
$$
Alors 
$$
\L z^*\big(\psi^k(\Lambda_{-1}(V))\big)=0.
$$
\label{implem}
\end{lemme}
\beginProof
Soit 
$$
z^*V\to\CO\to 0
$$
la restriction de la suite exacte à $Z$. Soit 
$W$ le noyau du morphisme $z^*V\to\CO$. Par construction, 
$W$ est localement libre. Par ailleurs, on l'équation
$$
\psi^k(\Lambda_{-1}(z^*V))=\psi^k(\Lambda_{-1}(W))\otimes
\psi^k(\Lambda_{-1}(\CO)).
$$
Comme $\psi^k(\Lambda_{-1}(\CO))=0$, on déduit le lemme.
\endProof
\begin{lemme}
Soit $i\not=j$ et soit $Z$ un sous-schéma fermé de $D_i\cap D_j$. Alors 
pour tout $n\in\mZ$, on a
$$
\theta^{k}(\CO_{Z})^{n}
\otimes\psi^k(\Lambda_{-1}(\Omega^1_{X/Y}(\log)))=
\psi^k(\Lambda_{-1}(\Omega^1_{X/Y}(\log)))
$$
dans $\K_0(X)[{1\over k}].$
\end{lemme}
\beginProof
Soit $z:Z\hookrightarrow X$ l'immersion 
canonique. Soit $j:U:=X\backslash Z\hookrightarrow X$ l'immersion ouverte 
du complément de $C$. 
On rappelle que l'on dispose d'une suite exacte
$$
\K_0(Z)\stackrel{\R z_*}{\to}\K_0(X)\stackrel{\L j^*}{\to}\K_0(U)\to 0.
$$
On calcule que 
$$
\L j^*\big(\theta^{k}(\CO_{Z})^{n}
-1\big)=0.
$$
On en déduit qu'il existe un élément $x\in\K_0(Z)[{1\over k}]$ tel que
$$
\R z_*(x)=\theta^{k}(\CO_{Z})^{n}
-1.
$$
Ainsi
$$
\big(\theta^{k}(\CO_{Z})^{n}
-1\big)\otimes\psi^k(\Lambda_{-1}(\Omega^1_{X/Y}(\log)))=
\R z_*(x\otimes\L z^*(\psi^k(\Lambda_{-1}(\Omega^1_{X/Y}(\log))))).
$$
Les Lemmes \ref{firstlem} et \ref{implem} impliquent alors que 
\begin{equation*}
\L z^*(\psi^k(\Lambda_{-1}(\Omega^1_{X/Y}(\log))))=0.
\end{equation*}
\endProof

L'équation \refeq{pivotadams} suit maintenant du lemme suivant.

\begin{lemme}
Dans $\K'_0(X)$, le faisceau cohérent 
${\cal O}_{D_0}/f^*{\cal O}_{E_0}$ est une combinaison linéaire 
de faisceaux structuraux $\CO_Z$ de sous-schémas fermés 
$Z$ de $\coprod_{i\not=j}({D_i\cap D_j})$. 
\label{comblem}
\end{lemme}
\beginProof
On commence par rappeler l'énoncé suivant. Soit 
$T$ un schéma réduit et noethérien. Alors 
$\K'_0(T)$ est engendré par les classes 
$\CO_Z$, où $Z$ est un sous-schéma fermé intègre de $T$. 
Pour démontrer cela, on suppose tout d'abord que l'énoncé est vérifié 
si l'on remplace $T$ par $T'$, où $T'$ est un sous-schéma fermé réduit tel que 
$T'\not=T$. On remarque que pour tout faisceau cohérent $F$ sur $T$ il existe 
un ouvert  non-vide $U\subseteq T$ tel que 
$F|_U$ est une somme directe de faisceaux $\CO_U$. Par ailleurs, on 
a une injection $\CO_T\hookrightarrow\oplus_{C\in\Irr(T)}\CO_{C}$ et 
le conoyau du morphisme $\CO_C\rightarrow\oplus_{C\in\Irr(T)}\CO_{C}$ est supporté 
par un sous-schéma fermé qui ne coïncide pas pas avec $T$. 
Le principe d'induction noethérienne permet alors de conclure la 
démonstration de l'énoncé.

Le lemme est une conséquence de ce dernier énoncé et du fait 
que le faisceau ${\cal O}_{D_0}/f^*{\cal O}_{E_0}$ est supporté sur 
le sous-schéma fermé $\cup_{i\not=j}{D_i\cap D_j}$.  
\endProof
Les équations \refeq{probtermp} et \refeq{pivotadams} mènent pour finir
à l'équation 
$$
\theta^k(f)^{-1}\psi^k(\Lambda_{-1}(\Omega^1_{X/Y}(\log)))=
\theta^k(\Omega^1_{X/Y}(\log))^{-1}\psi^k(\Lambda_{-1}(\Omega^1_{X/Y}(\log)))
$$
qui conclut notre analyse du terme \refeq{probterm}. 
On insère cette dernière égalité dans l'équation \refeq{eqadams}. On 
obtient 
\begin{eqnarray*}
\psi^k(\R f_*(\Lambda_{-1}(\Omega^1_{X/Y}(\log)))&=&
\R f_*(\theta^k(f)^{-1}\psi^k(\Lambda_{-1}(\Omega^1_{X/Y}(\log))))\\
&=&
\R f_*(\theta^k(\Omega^1_{X/Y}(\log))^{-1}\psi^k(\Lambda_{-1}(\Omega^1_{X/Y}(\log))))
=
\R f_*(\Lambda_{-1}(\Omega^1_{X/Y}(\log)))
\end{eqnarray*}
où on a utilisé l'équation 
$$
\psi^k(V)(\Lambda_{-1}(V))=\theta^k(V)\Lambda_{-1}(V)
$$
valable dans $\K_0(X)$ pour tout faisceau localement libre $V$ sur $X$ 
(cf. \cite[Prop. 7.3, p. 269]{Atiyah-Tall-Group}). On obtient ainsi que  
$$
\rc_{\rm tot}\big(\psi^k(\R f_*(\Lambda_{-1}(\Omega^1_{X/Y}(\log))))=
\rc_{\rm tot}\big(\psi^k(\R f_*(\Lambda_{-1}(\Omega^1_{X/Y}(\log))))
=\rc_{\rm tot}(\R f_*(\Lambda_{-1}(\Omega^1_{X/Y}(\log))))
$$
L'existence de la suite spectrale de Hodge vers de Rham montre 
que 
$$
\R f_*(\Lambda_{-1}(\Omega^1_{X/Y}(\log)))=\sum_{k\geqslant 0}(-1)^k\H^k_\dR(X/Y)(\log)
$$
et on conclut que  
$$
(k^t-1)\cdot\rc_{t}(\sum_{k\geqslant 0}(-1)^k\H^k_\dR(X/Y)(\log))=0.
$$
dans ${\rm CH}^t(Y)[{1\over k}]$, pour tout $t\geqslant 0$. Comme 
ce résultat est valable pour tout $k\geqslant 2$, 
le Théorème \ref{mainth} (b) est maintenant un 
conséquence du Lemme \ref{adamslem}.

{\it Passons maintenant à la démonstration de (e).}

{\it On démontre d'abord la première assertion.  }
Soit $Y':=\Spec\ f_*\CO_X$ et soit 
$g:Y'\to Y$ le morphisme naturel.
Comme le faisceau $f_*\CO_X$ est 
localement libre (voir Théorème \ref{ill-deg}), le morphisme 
$g$ est plat et fini.
 Soit $\mtr{y}\to Y$ un point géométrique. Comme les faisceaux $\R^i f_*\CO_X$ sont 
aussi localement libres (toujours le Théorème \ref{ill-deg}), 
le théorème de semi-continuité (cf. par ex. \cite[III, th. 12.11]{Hartshorne-Algebraic}) montre que 
$$
Y'_{\mtr{y}}\simeq\Spec\ \H^0(X_{\mtr{y}},\CO_{X_{\mtr{y}}}).
$$
Comme $X_{\mtr{y}}$ est réduit (voir sous-section \ref{morss}), on voit 
que les fibres géométriques de $g$ sont réduites. 
On déduit que $g$ est étale (cf. \cite[IV, par. 17, Cor. 17.6.2, c'')]{EGA}). 
On remarque maintenant que par construction, on a un isomorphisme
$$
\H^0_\dR(X/Y)(\log)\simeq g_*\CO_{Y'}.
$$
Pour démontrer la première assertion de (e), on est donc ramené à la 
démontrer dans le cas où $f$ est un morphisme étale et fini. Dans cas cas-là, 
l'assertion résulte de (d).

{\it Passons à la deuxième assertion.}

Il est montré dans \cite[1.3, p. 144]{Illusie-Reduction} que si les fibres 
de $f$ sont de dimension $1$, 
alors le faisceau 
dualisant relatif de $f$ est canoniquement isomorphe 
à $\Omega^1_{X/Y}(\log)$. Par dualité 
de Grothendieck, on obtient un isomorphisme 
canonique
$$
(\R^1 f_*\Omega^1_{X/Y}(\log))^\vee\simeq\R^0 f_*\CO_X.
$$
Par ailleurs, comme la suite spectrale de Hodge vers de Rham logarithmique
dégénère, on a des isomorphismes
$$
\R^0f_*\CO_X\simeq\H^0_\dR(X/Y)(\log)
$$
et 
$$
\R^1 f_*\Omega^1_{X/Y}(\log)\simeq\H^2_\dR(X/Y)(\log).
$$
On applique maintenant (d). On obtient
$$
\N_t\cdot\rc_t\big(-\H^1_\dR(X/Y)(\log)+\R^0 f_*\CO_X+(\R^0 f_*\CO_X)^\vee\big)=0.
$$
La deuxième assertion résulte maintenant de la première assertion et du lemme suivant.

\begin{lemme}
Soit $e,e'\in\K_0(Y)$. On suppose que 
$\N_l\cdot\rc_l(e)=0$ et $\N_l\cdot\rc_l(e')=0$ pour tout 
$l\geqslant 0$. Alors $\N_l\cdot\rc_l(e+e')=\N_l\cdot\rc_l(-e)=0$ 
pour tout $l\geqslant 0$. 
\end{lemme}
\beginProof
Nous allons utiliser la divisibilité 
$$
\pgcd(\N_{j},\N_{k})|\N_{j+k}
$$
valable pour tout $j,k\geqslant 0$ (voir \cite[Rem. 19.5, p. 66]{Fulton-MacPherson-Char}). 
On calcule
$$
\N_l\cdot\rc_{l}(e+e')=\sum_{j+k=l}(\N_l\cdot \rc_j(e)\rc_k(e'))
$$
et comme par hypothèse $\pgcd(\N_j,\N_k)\cdot \rc_j(e)\rc_k(e')=0$ pour tout $j,k\geqslant 0$, on conclut que $\N_l\cdot\rc_{l}(e+e')=0$. La démonstration du fait que 
$\N_l\cdot\rc_l(-e)=0$ est similaire. 
\endProof

\subsection{Démonstration de (f)}

Soit $G$ un groupe. Pour tout anneau commutatif $A$, nous 
noterons $\R_A(G)$ le groupe de Grothendieck des 
représentations linéaires de $G$ dans des $A$-modules libres. 
On rappelle qu'une représentation linéaire de $G$ dans 
un $A$-module est la donnée d'un $A$-module $M$, d'un nombre 
entier $n\geqslant 0$ et 
d'un homomorphisme de groupes $G\to\GL_n(A)$. 
Le produit tensoriel des représentations munit le groupe $\R_A(G)$ 
d'une structure d'anneau commutatif et le produit extérieur 
le munit d'une structure de $\lambda$-anneau (voir par ex. \cite{Serre-Groupes} pour tout cela).
  
Soit $M$ une variété quasi-projective lisse sur $\mC$. 
On munit $M(\mC)$ de sa structure canonique d'espace analytique 
complexe. Si $V$ est un faisceau cohérent localement libre sur $M$, 
on écrira $V^\an$ pour le fibré vectoriel holomorphe associé sur 
$M(\mC)$; si $\nabla$ est une connexion sur 
$V$, on écrira $\nabla^\an$ pour la connexion analytique 
sur $V^\an$ associée à $\nabla$. 
On rappelle que Deligne a démontré le théorème suivant: le foncteur 
$(V,\nabla)\mapsto(V^\an,\nabla^\an)$ que l'on vient de décrire 
induit une équivalence entre la catégorie additive 
des faisceaux cohérents localement libres munis d'un 
connexion intégrable régulière et la 
catégorie additive des fibrés vectoriels holomorphes 
munis d'un connexion (analytique) intégrable (cf. \cite[II,Th. 5.9, p. 97]{Deligne-Equations}). 
Pour la définition d'une connexion intégrable régulière, 
voir \cite[II, Prop. 4.4]{Deligne-Equations} (cette notion est due 
à Griffiths). Deligne démontre aussi 
que la catégorie  des faisceaux cohérents localement libres munis d'un 
connexion intégrable régulière est 
fermée par produit tensoriel et puissances extérieures 
(cf. \cite[II, Prop. 4.6, p. 90]{Deligne-Equations}).

Soit maintenant $m\in M(\mC)$ un point de base. 
Comme la catégorie des fibrés vectoriels holomorphes  
sur $M(\mC)$ munis d'une connexion intégrable est 
équivalente à la catégorie des représentations 
linéaires du groupe fondamental $\pi_1(M(\mC),m))$, 
le théorème de Deligne fournit un morphisme canonique de groupes
$$
\R_{\mC}(\pi_1(M(\mC),m))\to\K_0(M). 
$$
Par construction, ce morphisme est un morphisme de 
$\lambda$-anneaux.  

\begin{prop}
On suppose que $G$ est fini. 
Si $(k,\#G)=1$, l'opération 
d'Adams $\psi^k:\R_\mQ(G)\to\R_\mQ(G)$ est l'identité.
\label{propQ}
\end{prop}
\beginProof
On suppose tout d'abord que $G$ est cyclique. 
Soit $n:=\#G$.  Soit 
$\rho:G\to\GL(V)$ une représentation de $G$ dans 
un espace vectoriel complexe $V$. 
Par construction, on a 
$\Tr(\rho(g))\in\mQ(\mu_n)$ pour tout $g\in G$. 
Pour tout $k\in(\mZ/n\mZ)^*$ notons $\sigma_k\in\Gal(\mQ(\mu_n)|\mQ)$ 
l'élément associé. Nous allons démontrer  que 
\begin{equation}
\Tr(\psi^k(\rho)(g))=\sigma_k(\Tr(\rho(g))).
\label{eqaux}
\end{equation}
Ici $\psi^k:\R_\mC(G)\to\R_\mC(G)$ est la $k$-ème 
opération d'Adams sur $\R_\mC(G)$. 
Comme les deux côtés de cette égalité sont 
additifs pour les sommes directes de représentations, 
on est ramené au cas où $\rho$ est de dimension $1$. 
Dans ce cas-là, l'égalité est une conséquence des définitions.

Revenons aux hypothèses de la proposition. 
Vu que la trace et $\psi^k(\cdot)$ sont additifs pour les 
sommes directes de représentations, il suffit, pour démontrer 
la proposition, de montrer 
l'égalité 
\begin{equation}
\Tr(\psi^k(\rho)(g))=\Tr(\rho(g))
\label{eqaux2}
\end{equation}
pour toute représentation linéaire $\rho:G\to\GL(V)$ dans un 
espace vectoriel $V$ sur $\mQ$. Pour un $g$ donné, 
les deux membres de cette égalité restent inchangés si 
l'on remplace $G$ par un de ses sous-groupes contenant $g$. 
On peut ainsi sans restriction de généralité supposer que 
$G$ est un groupe cyclique engendré par $g$. 
Comme le morphisme naturel d'extension des scalaires 
$\R_\mQ(G)\to\R_\mC(G)$ préserve les traces et est un morphisme de 
$\lambda$-anneaux, il suffit de démontrer l'équation 
$$
\Tr(\psi^k(\rho_\mC)(g))=\Tr(\rho_\mC(g))
$$
pour la représentation complexe $\rho_\mC$ associée à $\rho$. 
Comme $\Tr(\rho_\mC(g))\in\mQ$, cette dernière équation 
est une conséquence de \refeq{eqaux}. 
\endProof

On peut maintenant démontrer (f).  Soit $y\in Y(\mC)$ un point de base. 
Soit $V:=(\R^j f(\mC)_*\mQ)_y$. Soit $G\subseteq\GL(V)$ l'image 
de la représentation de monodromie $\pi_1(Y(\mC),y)\to \GL(V)$
associée à $\R^j f(\mC)_*\mQ$. Le fibré vectoriel holomorphe 
associé à cette représentation est isomorphe 
au fibré $\H^j_\dR(X/Y)^\an$. Par ailleurs, il existe 
une connexion intégrable régulière $\nabla$ sur 
$\H^j_\dR(X/Y)$ telle que $\nabla^\an$ est isomorphe à la connexion 
intégrable (analytique) sur $\H^j_\dR(X/Y)^\an$ qui est induite par 
la représentation de monodromie.  
Ceci est une conséquence de résultats de Grifiths (cf. \cite{Griffiths-Periods-III}) et Katz-Oda 
(cf. \cite{Katz-Oda-Differentiation}). Voir \cite[III, Th. 7.9]{Deligne-Equations} pour la 
démonstration. On appelle {\it connexion de 
Gauss-Manin} la connexion $\nabla$.

Soit $\rho:G\to\GL(V)$ 
le morphisme d'inclusion. On dispose de morphismes naturels 
de $\lambda$-anneaux
$$
\R_\mQ(G)\to\R_{\mQ}(\pi_1(Y(\mC),y))\to\R_{\mC}(\pi_1(Y(\mC),y))\to\K_0(Y).
$$
Soit $r:\R_\mQ(G)\to\K_0(Y)$ le morphisme composé. 
Les propriétés de la connexion de Gauss-Manin mentionnées plus haut impliquent que 
$r(\rho)=\H^j_{\dR}(X/Y)$ dans $\K_0(Y)$. 
En tenant compte de la Proposition \ref{propQ}, on voit ainsi que 
$$
0=\rc_t(r(\psi^p(\rho)-\rho))=\rc_t(\psi^p(r(\rho))-r(\rho))=(p^t-1)\rc_t(r(\rho))=
(p^t-1)\rc_t(\H^j_\dR(X/Y))
$$
pour presque tout nombre premier $p$. On conclut au moyen du Lemme
\ref{adamslem}. 

\section{Démonstration de (g) \& (h)}

{\it Commen\c{c}ons par la démonstration de (g).} 
On peut supposer, sans restriction de généralité, que $Y$ est connexe. 
Soit encore une fois $L$ un corps de type fini sur $\mQ$ (comme corps) tel que 
le triplet $f:X\to Y$ possède un modèle $f_0:X_0\to Y_0$ sur $L$.  
Soit 
$T$ un schéma affine, intègre, lisse et de type fini sur  
$\mZ$, dont le corps de fonctions est isomorphe à $L$. Comme avant, quitte 
à remplacer $T$ par l'un de ses ouverts, on peut supposer qu'il existe 
des modèles lisses et projectifs 
$\widetilde{X}$ et $\widetilde{Y}$ de $X_0$ et $Y_0$ sur $T$ et 
qu'il existe un $T$-morphisme lisse et projectif
$\widetilde{f}:\widetilde{X}\ra\widetilde{Y}$ qui est 
un modèle de $f_0$. 
Toujours quitte à réduire la taille de $T$, on peut supposer que la suite spectrale de Hodge vers de Rham 
de $\widetilde{f}:\widetilde{X}\ra\widetilde{Y}$ dégénère et que les faisceaux de cohomologie 
de de Rham relative correspondants sont localement libres. Vu les hypothèeses 
de (g), on peut aussi, quitte à restreindre $T$, supposer 
que la classe de Chern $\rc_t(\H^j_\dR(\wt{X}/\wt{Y})(\log))\in\CHOW^t(\wt{Y})$ est de torsion. 

Soit ${\mfp}$ un point de $T$ et soit $\mtr{\mfp}\to T$ le point géométrique 
associé à une clôture algébrique de $\kappa(\mfp)$. Soit $\mtr{L}$ un clôture algébrique de $L$. 
Soit $l$ un 
nombre premier différent de la caractéristique résiduelle de $\mfp$. 
Nous allons décrire la construction d'un {\it homomorphisme de 
 spécialisation} $$\sigma:\CHOW^t(Y_{0,\mtr{L}})[l^\infty]\to \CHOW^t(\widetilde{Y}_{\mtr{\mfp}})[l^\infty].$$ Cet homomorphisme n'est pas canonique et dépend du choix 
 de certains plongements (voir aussi \cite[Introduction]{Schoen-Spec} pour une autre description de 
 $\sigma$).
 Soit $T_1\to T$ l'éclatement de $T$ en $\mfp$. Soit 
 $\mfp_1$ le point générique de la fibre spéciale 
 de $T_1$. Soit $R$ l'anneau local en $\mfp_1$ et soit $J$ l'idéal maximal 
 de $R$. Cet anneau est par construction un anneau de 
 valuation discret et $\Frac(R)=L$. Soit enfin $\wt{Y}_1:=\wt{Y}\times_T\Spec\ R$. 
 On dispose d'un homomorphisme  de spécialisation
 $$
 \sigma_L:\CHOW^t(Y_0)\to \CHOW^t(\widetilde{Y}_{1,\kappa(\mfp_1)})
$$
 obtenu en associant à chaque sous-schéma fermé intègre 
 $Z$  
 de $Y_0$ la classe de la restriction de $\Zar_{\wt{Y}_1}(Z)$ 
 à $\widetilde{Y}_{1,\kappa(\mfp_1)}$. Voir \cite[20.3]{Fulton-Intersection} pour plus de détails. 
 Soit $\mtr{J}$ un idéal maximal de la clôture entière de $R$ dans $\mtr{L}$. 
 Pour chaque extension $M$ de $L$ plongée dans $\mtr{L}$, soit 
 $R'_M$ la clôture intégrale de $R$ dans $M$ et soit $R_M$ la localisation 
 de $R'_M$ en $\mtr{J}\cap R'_M$. Notons $J_M$ l'idéal maximal de $R_M$ 
 et $\lambda_M$ le corps résiduel associé. 
 Si $M$ est un extension finie de $L$ alors l'anneau $R_M$ est à nouveau un anneau 
 de valuation discret. On dispose ainsi pour chaque $M$ fini sur $L$ 
 d'un homomorphisme de spécialisation 
 $\sigma_M:\CHOW^t(Y_{0,M})\to \CHOW^t(\widetilde{Y}_{1,\lambda_M})$ 
 et si $M'\supseteq M$ est une paire emboîtée d'extensions, on 
 a un diagramme commutatif naturel
 $$
 \xymatrix{
 \CHOW^t(Y_{0,M})\ar[r]^{\sigma_M}\ar[d] & \CHOW^t(\widetilde{Y}_{1,\lambda_M})\ar[d]\\
  \CHOW^t(Y_{0,M'})\ar[r]^{\sigma_{M'}}& \CHOW^t(\widetilde{Y}_{1,\lambda_{M'}})
  }
  $$
On vérifie que les homomorphismes naturels 
$$\indlim_{M, [M:L]<\infty}\CHOW^t(Y_{0,M})\to\CHOW^t(Y_{0,\mtr{L}})$$ 
et $$\indlim_{M, [M:L]<\infty}\CHOW^t(\wt{Y}_{1,\lambda_M})\to\CHOW^t(\wt{Y}_{1,\lambda_{\mtr{L}}}).$$ 
sont des isomorphismes. On obtient ainsi un homomorphisme 
de spécialisation
$$
\sigma_{\mtr{L}}:\CHOW^t(Y_{0,\mtr{L}})\to \CHOW^t(\wt{Y}_{1,\lambda_{\mtr{L}}})
$$
Choisissons un plongement $\kappa(\mtr{\mfp})\hookrightarrow \lambda_{\mtr{L}}$ compatible au plongement canonique $\kappa(\mfp)\hookrightarrow \kappa(\mfp_1)$. 
On peut montrer que l'homomorphisme canonique
$$
\CHOW^t(\wt{Y}_{\mtr{\mfp}})\to\CHOW^t(\wt{Y}_{1,\lambda_{\mtr{L}}})
$$
induit un isomorphisme 
$$
\tau_1:\CHOW^t(\wt{Y}_{\mtr{\mfp}})[l^\infty]\to\CHOW^t(\wt{Y}_{1,\lambda_{\mtr{L}}})[l^\infty].
$$
De même, le choix d'un $L$-plongement $\mtr{L}\hookrightarrow K$ 
induit un isomorphisme
$$
\tau_2:\CHOW^t(Y_{0,\mtr{L}})[l^\infty]\to\CHOW^t(Y)[l^\infty].
$$
Pour ce résultat non-trivial, voir \cite{Lecomte-Rig}. 
On définit pour finir 
$$
\sigma:=(\tau_1)^{-1}\circ\sigma_{\mtr{L}}\circ(\tau_2)^{-1}.
$$
Par ailleurs, le théorème 
du changement de base propre en cohomologie étale 
implique l'existence d'un homomorphisme de 
spécialisation canonique 
$$
\sigma_\et:\H^{2t-1}_\et(Y,\mQ_l/\mZ_l(t))\to
\H^{2t-1}_\et(\widetilde{Y}_{\mtr{\mfp}},\mQ_l/\mZ_l(t))
$$
qui est un isomorphisme. Le 
 théorème de spécialisation 
pour l'application d'Abel-Jacobi de Bloch (cf. \cite[Prop. 3.8]{Bloch-Torsion}) implique 
alors que $\lambda^t_l\circ \sigma=\sigma_\et\circ\lambda^t_l$. 

Si $G$ est un groupe abélien et $g\in G$ est un élément 
de torsion, notons $g[l]$ la partie $l$-primaire de $g$. 

\begin{lemme}
Soit $F$  un fibré cohérent localement libre sur $\wt{Y}$. Supposons 
que $\rc_t(F)$ est de torsion. Alors  
$$
\sigma(\rc_t(F_{K})[l])=\rc_t(F_{\mtr{\mfp}})[l]
$$
\label{reader-lem}
\end{lemme}
Le Lemme \ref{reader-lem} est une conséquence du fait 
que les classes de Chern commutent aux images réciproques de fibrés localement 
libres par des morphismes de schémas. Les détails de la démonstration sont 
laissés au lecteur. 

En utilisant le Lemme \ref{reader-lem}, 
on peut, par un argument tout semblable à celui apparaissant dans la 
preuve de (a), conclure que pour tout nombre premier $l$, l'égalité
$$
(1-p^t)\cdot\lambda^t_l(\rc^t(\H^j_\dR(X/Y)(\log))[l])=0
$$
est vérifiée pour presque tout nombre premier $p$. On invoque 
alors le Lemme \ref{adamslem} et on voit que 
$$
\N_t\cdot\lambda^t_l(\rc^t(\H^j_\dR(X/Y)(\log))[l])=
\N_t\cdot\lambda^t_l(\rc^t(\H^j_\dR(X/Y)(\log)))[l]=0
$$
pour tout nombre premier $l$. Autrement dit,
$$
\N_t\cdot\lambda^t(\rc^t(\H^j_\dR(X/Y)(\log)))=0, 
$$
ce qu'on voulait démontrer.

{\it Passons maintenant à la démonstration de (h)}. Dans 
\cite[Par. 1.2, Cor. 4]{Colliot-Sansuc-Torsion}, il est démontré que l'application
$$
\lambda^2:\Tor(\CHOW^2(Y))\to\bigoplus_{l\ {\rm premier}}\H^{3}_\et(Y,\mQ_l/\mZ_l(2)).
$$
est injective. On voit donc que (g) implique (h).


\begin{bibdiv}
\begin{biblist}
\bib{Adams-J(X)}{misc}{
  author={Adams, J. F.},
  title={On the groups $J(X)$. II, III, IV. {\rm Topology {\bf 3 \& 5}, 1965-1966}},
}

\bib{Atiyah-Tall-Group}{article}{
  author={Atiyah, M. F.},
  author={Tall, D. O.},
  title={Group representations, $\lambda $-rings and the $J$-homomorphism},
  journal={Topology},
  volume={8},
  date={1969},
  pages={253--297},
  issn={0040-9383},
}

\bib{BFQ}{article}{
  author={Baum, Paul},
  author={Fulton, William},
  author={Quart, George},
  title={Lefschetz-Riemann-Roch for singular varieties},
  journal={Acta Math.},
  volume={143},
  date={1979},
  number={3-4},
  pages={193--211},
  issn={0001-5962},
}

\bib{Bismut-Eta}{article}{
  author={Bismut, J. M.},
  title={Eta invariants, differential characters and flat vector bundles},
  note={With an appendix by K. Corlette and H. Esnault},
  journal={Chinese Ann. Math. Ser. B},
  volume={26},
  date={2005},
  number={1},
  pages={15--44},
  issn={0252-9599},
}

\bib{Bloch-Torsion}{article}{
  author={Bloch, S.},
  title={Torsion algebraic cycles and a theorem of Roitman},
  journal={Compositio Math.},
  volume={39},
  date={1979},
  number={1},
  pages={107--127},
  issn={0010-437X},
}

\bib{Bloch-Esnault-Algebraic}{article}{
  author={Bloch, Spencer},
  author={Esnault, H{\'e}l{\`e}ne},
  title={Algebraic Chern-Simons theory},
  journal={Amer. J. Math.},
  volume={119},
  date={1997},
  number={4},
  pages={903--952},
  issn={0002-9327},
}

\bib{Colliot-Sansuc-Torsion}{article}{
  author={Colliot-Th{\'e}l{\`e}ne, Jean-Louis},
  author={Sansuc, Jean-Jacques},
  author={Soul{\'e}, Christophe},
  title={Torsion dans le groupe de Chow de codimension deux},
  language={French},
  journal={Duke Math. J.},
  volume={50},
  date={1983},
  number={3},
  pages={763--801},
  issn={0012-7094},
}

\bib{Deligne-Equations}{book}{
  author={Deligne, Pierre},
  title={\'Equations diff\'erentielles \`a points singuliers r\'eguliers},
  language={French},
  series={Lecture Notes in Mathematics, Vol. 163},
  publisher={Springer-Verlag},
  place={Berlin},
  date={1970},
  pages={iii+133},
}

\bib{Demazure-Cours}{book}{
  author={Demazure, Michel},
  title={Cours d'alg\`ebre},
  language={French, with French summary},
  series={Nouvelle Biblioth\`eque Math\'ematique [New Mathematics Library], 1},
  note={Primalit\'e. Divisibilit\'e. Codes. [Primality. Divisibility. Codes]},
  publisher={Cassini, Paris},
  date={1997},
  pages={xviii+302},
  isbn={2-84225-000-1},
}

\bib{Ekedahl-VDG}{article}{
  author={Ekedahl, Torsten},
  author={van der Geer, Gerard},
  title={The order of the top Chern class of the Hodge bundle on the moduli space of abelian varieties},
  journal={Acta Math.},
  volume={192},
  date={2004},
  number={1},
  pages={95--109},
  issn={0001-5962},
}

\bib{Esnault-Viehweg-Chern}{article}{
  author={Esnault, H{\'e}l{\`e}ne},
  author={Viehweg, Eckhart},
  title={Chern classes of Gauss-Manin bundles of weight 1 vanish},
  journal={$K$-Theory},
  volume={26},
  date={2002},
  number={3},
  pages={287--305},
  issn={0920-3036},
}

\bib{Evens-Kahn-Chern}{article}{
  author={Evens, Leonard},
  author={Kahn, Daniel S.},
  title={Chern classes of certain representations of symmetric groups},
  journal={Trans. Amer. Math. Soc.},
  volume={245},
  date={1978},
  pages={309--330},
  issn={0002-9947},
}

\bib{Fulton-Intersection}{book}{
  author={Fulton, William},
  title={Intersection theory},
  series={Ergebnisse der Mathematik und ihrer Grenzgebiete. 3. Folge. A Series of Modern Surveys in Mathematics [Results in Mathematics and Related Areas. 3rd Series. A Series of Modern Surveys in Mathematics]},
  volume={2},
  edition={2},
  publisher={Springer-Verlag},
  place={Berlin},
  date={1998},
  pages={xiv+470},
  isbn={3-540-62046-X},
  isbn={0-387-98549-2},
}

\bib{Fulton-MacPherson-Char}{article}{
  author={Fulton, William},
  author={MacPherson, Robert},
  title={Characteristic classes of direct image bundles for covering maps},
  journal={Ann. of Math. (2)},
  volume={125},
  date={1987},
  number={1},
  pages={1--92},
  issn={0003-486X},
}

\bib{Griffiths-Periods-III}{article}{
  author={Griffiths, Phillip A.},
  title={Periods of integrals on algebraic manifolds. III. Some global differential-geometric properties of the period mapping},
  journal={Inst. Hautes \'Etudes Sci. Publ. Math.},
  number={38},
  date={1970},
  pages={125--180},
  issn={0073-8301},
}

\bib{FGA-232}{misc}{
  author={Grothendieck, A.},
  title={Technique de descente et Th\'eor\`emes d'existence en g\'eom\'etrie alg\'ebrique. V. Les sch\'emas de Picard : th\'eor\`emes d'existence. {\rm S\'eminaire Bourbaki, 7 (1961-1962), Expos\'e No. 232.}},
}

\bib{Grothendieck-Classes}{article}{
  author={Grothendieck, A.},
  title={Classes de Chern et repr\'esentations lin\'eaires des groupes discrets},
  language={French},
  conference={ title={Dix Expos\'es sur la Cohomologie des Sch\'emas}, },
  book={ publisher={North-Holland}, place={Amsterdam}, },
  date={1968},
  pages={215--305},
}

\bib{EGA}{article}{
  author={Grothendieck, A.},
  status={{\it \'El\'ements de g\'eom\'etrie alg\'ebrique.} { Inst. Hautes \'Etudes Sci. Publ. Math.} {\bf 4, 8, 11, 17, 20, 24, 28, 32} (1960-1967).},
}

\bib{SGA6}{book}{
  title={Th\'eorie des intersections et th\'eor\`eme de Riemann-Roch},
  note={S\'eminaire de G\'eom\'etrie Alg\'ebrique du Bois-Marie 1966--1967 (SGA 6); Dirig\'e par P. Berthelot, A. Grothendieck et L. Illusie. Avec la collaboration de D. Ferrand, J. P. Jouanolou, O. Jussila, S. Kleiman, M. Raynaud et J. P. Serre; Lecture Notes in Mathematics, Vol. 225},
  publisher={Springer-Verlag},
  place={Berlin},
  date={1971},
  pages={xii+700},
}

\bib{SGA3-1}{book}{
  title={Sch\'emas en groupes. I: Propri\'et\'es g\'en\'erales des sch\'emas en groupes},
  language={French},
  series={S\'eminaire de G\'eom\'etrie Alg\'ebrique du Bois Marie 1962/64 (SGA 3). Dirig\'e par M. Demazure et A. Grothendieck. Lecture Notes in Mathematics, Vol. 151},
  publisher={Springer-Verlag},
  place={Berlin},
  date={1970},
  pages={xv+564},
}

\bib{SGA4/2}{book}{
  author={Deligne, P.},
  title={Cohomologie \'etale},
  series={Lecture Notes in Mathematics, Vol. 569},
  note={S\'eminaire de G\'eom\'etrie Alg\'ebrique du Bois-Marie SGA 4${1\over 2}$; Avec la collaboration de J. F. Boutot, A. Grothendieck, L. Illusie et J. L. Verdier},
  publisher={Springer-Verlag},
  place={Berlin},
  date={1977},
  pages={iv+312pp},
}

\bib{SGA4.3}{book}{
  title={Th\'eorie des topos et cohomologie \'etale des sch\'emas. Tome 3},
  language={French},
  series={Lecture Notes in Mathematics, Vol. 305},
  note={S\'eminaire de G\'eom\'etrie Alg\'ebrique du Bois-Marie 1963--1964 (SGA 4); Dirig\'e par M. Artin, A. Grothendieck et J. L. Verdier. Avec la collaboration de P. Deligne et B. Saint-Donat},
  publisher={Springer-Verlag},
  place={Berlin},
  date={1973},
  pages={vi+640},
}

\bib{SGA5}{book}{
  title={Cohomologie $l$-adique et fonctions $L$},
  language={French},
  series={Lecture Notes in Mathematics, Vol. 589},
  note={S\'eminaire de G\'eometrie Alg\'ebrique du Bois-Marie 1965--1966 (SGA 5); Edit\'e par Luc Illusie},
  publisher={Springer-Verlag},
  place={Berlin},
  date={1977},
  pages={xii+484},
  isbn={3-540-08248-4},
}

\bib{Hartshorne-Algebraic}{book}{
  author={Hartshorne, Robin},
  title={Algebraic geometry},
  note={Graduate Texts in Mathematics, No. 52},
  publisher={Springer-Verlag},
  place={New York},
  date={1977},
  pages={xvi+496},
  isbn={0-387-90244-9},
}

\bib{Illusie-Reduction}{article}{
  author={Illusie, Luc},
  title={R\'eduction semi-stable et d\'ecomposition de complexes de de Rham \`a\ coefficients},
  language={French},
  journal={Duke Math. J.},
  volume={60},
  date={1990},
  number={1},
  pages={139--185},
  issn={0012-7094},
}

\bib{Katz-Oda-Differentiation}{article}{
  author={Katz, Nicholas M.},
  author={Oda, Tadao},
  title={On the differentiation of de Rham cohomology classes with respect to parameters},
  journal={J. Math. Kyoto Univ.},
  volume={8},
  date={1968},
  pages={199--213},
  issn={0023-608X},
}

\bib{Koeck-Grothendieck}{article}{
  author={K{\"o}ck, Bernhard},
  title={The Grothendieck-Riemann-Roch theorem for group scheme actions},
  language={English, with English and French summaries},
  journal={Ann. Sci. \'Ecole Norm. Sup. (4)},
  volume={31},
  date={1998},
  number={3},
  pages={415--458},
  issn={0012-9593},
}

\bib{Lecomte-Rig}{article}{
  author={Lecomte, Florence},
  title={Rigidit\'e des groupes de Chow},
  language={French},
  journal={Duke Math. J.},
  volume={53},
  date={1986},
  number={2},
  pages={405--426},
  issn={0012-7094},
}

\bib{Liu-Algebraic}{book}{
  author={Liu, Qing},
  title={Algebraic geometry and arithmetic curves},
  series={Oxford Graduate Texts in Mathematics},
  volume={6},
  note={Translated from the French by Reinie Ern\'e; Oxford Science Publications},
  publisher={Oxford University Press},
  place={Oxford},
  date={2002},
  pages={xvi+576},
  isbn={0-19-850284-2},
}

\bib{Maillot-Rossler-Conjectures}{article}{
  author={Maillot, Vincent},
  author={Roessler, Damien},
  title={Conjectures sur les d\'eriv\'ees logarithmiques des fonctions $L$ d'Artin aux entiers n\'egatifs},
  language={French, with English and French summaries},
  journal={Math. Res. Lett.},
  volume={9},
  date={2002},
  number={5-6},
  pages={715--724},
  issn={1073-2780},
}

\bib{Maillot-Rossler-On-the}{article}{
  author={Maillot, Vincent},
  author={Roessler, Damian},
  title={On the periods of motives with complex multiplication and a conjecture of Gross-Deligne},
  journal={Ann. of Math. (2)},
  volume={160},
  date={2004},
  number={2},
  pages={727--754},
  issn={0003-486X},
}

\bib{Milnor-Stasheff}{book}{
  author={Milnor, John W.},
  author={Stasheff, James D.},
  title={Characteristic classes},
  note={Annals of Mathematics Studies, No. 76},
  publisher={Princeton University Press},
  place={Princeton, N. J.},
  date={1974},
  pages={vii+331},
}

\bib{Pappas-Integral}{article}{
  author={Pappas, Georgios},
  title={Integral Grothendieck-Riemann-Roch theorem},
  journal={Invent. Math.},
  volume={170},
  date={2007},
  number={3},
  pages={455--481},
  issn={0020-9910},
}

\bib{Pink-Order}{article}{
  author={Pink, Richard},
  title={On the order of the reduction of a point on an abelian variety},
  journal={Math. Ann.},
  volume={330},
  date={2004},
  number={2},
  pages={275--291},
  issn={0025-5831},
}

\bib{Quillen-Higher}{article}{
  author={Quillen, Daniel},
  title={Higher algebraic $K$-theory. I},
  conference={ title={Algebraic $K$-theory, I: Higher $K$-theories (Proc. Conf., Battelle Memorial Inst., Seattle, Wash., 1972)}, },
  book={ publisher={Springer}, place={Berlin}, },
  date={1973},
  pages={85--147. Lecture Notes in Math., Vol. 341},
}

\bib{Reznikov-All}{article}{
  author={Reznikov, Alexander},
  title={All regulators of flat bundles are torsion},
  journal={Ann. of Math. (2)},
  volume={141},
  date={1995},
  number={2},
  pages={373--386},
  issn={0003-486X},
}

\bib{Rojtman-torsion}{article}{
  author={Rojtman, A. A.},
  title={The torsion of the group of $0$-cycles modulo rational equivalence},
  journal={Ann. of Math. (2)},
  volume={111},
  date={1980},
  number={3},
  pages={553--569},
  issn={0003-486X},
}

\bib{Rossler-Adams}{article}{
  author={Roessler, Damian},
  title={An Adams-Riemann-Roch theorem in Arakelov geometry},
  journal={Duke Math. J.},
  volume={96},
  date={1999},
  number={1},
  pages={61--126},
  issn={0012-7094},
}

\bib{Schoen-Spec}{article}{
  author={Schoen, Chad},
  title={Specialization of the torsion subgroup of the Chow group},
  journal={Math. Z.},
  volume={252},
  date={2006},
  number={1},
  pages={11--17},
  issn={0025-5874},
}

\bib{Serre-Groupes}{article}{
  author={Serre, Jean-Pierre},
  title={Groupes de Grothendieck des sch\'emas en groupes r\'eductifs d\'eploy\'es},
  language={French},
  journal={Inst. Hautes \'Etudes Sci. Publ. Math.},
  number={34},
  date={1968},
  pages={37--52},
  issn={0073-8301},
}

\bib{VDG-Cycles}{article}{
  author={van der Geer, Gerard},
  title={Cycles on the moduli space of abelian varieties},
  conference={ title={Moduli of curves and abelian varieties}, },
  book={ series={Aspects Math., E33}, publisher={Vieweg}, place={Braunschweig}, },
  date={1999},
  pages={65--89},
}

\end{biblist}
\end{bibdiv}
\end{document}